\documentclass[final]{autart}

\usepackage{url}
\usepackage{amsmath}
\usepackage{amssymb}
\usepackage{makeidx}
\usepackage{multicol}
\usepackage[bottom]{footmisc}
\usepackage{subfig}
\usepackage{graphicx}
\usepackage{color}
\usepackage{natbib}
\usepackage{epstopdf}
\usepackage{algorithm}
\usepackage{algcompatible}
\usepackage{calc}
\newtheorem{approxi}{Approximation}
\newtheorem{assump}{Assumption}
\newcommand{\proof}{\noindent {\bf Proof. }}
\newcommand{\app}{\begin{approxi}}
\newcommand{\eapp}{\end{approxi}}
\newcommand{\ass}{\begin{assump}}
\newcommand{\eass}{\end{assump}}
\newcommand{\teo}{\begin{thm}}
\newcommand{\eteo}{\end{thm}}
\newcommand{\corr}{\begin{cor}}
\newcommand{\ecorr}{\end{cor}}
\newcommand{\pro}{\begin{prop}}
\newcommand{\epro}{\end{prop}}
\newcommand{\lemma}{\begin{lem}}
\newcommand{\elemma}{\end{lem}}
\newcommand{\pb}{\begin{prob}}
\newcommand{\epb}{\end{prob}}
\newcommand{\df}{\begin{defn}}
\newcommand{\edf}{\end{defn}}
\newcommand{\rema}{\begin{rem}}
\newcommand{\erema}{\end{rem}}

\newcommand{\al}[1]{\begin{align} #1 \end{align}}
\newcommand{\nn}{\nonumber}

\newcommand{\mz}{\color{black}}

\newcommand{\tr}{\mathop{\rm tr}}  %traccia

\newcommand{\Cc}{ \mathcal{C}}

\newcommand{\Ec}{ \mathcal{E}}

\newcommand{\Gc}{ \mathcal{G}}

\newcommand{\Mc}{ \mathcal{M}}
\newcommand{\Nc}{ \mathcal{N}}

\newcommand{\Qc}{ \mathcal{Q}}

\newcommand{\Tc}{ \mathcal{T}}

\newcommand{\Vc}{ \mathcal{V}}

\newcommand{\Es}{ \mathbb{E}}

\newcommand{\Ns}{ \mathbb{N}}

\newcommand{\Rs}{ \mathbb{R}}

\newcommand{\Zs}{ \mathbb{Z}}

\newcommand{\Td}  {\mathrm{T}}

\begin{document}

\begin{frontmatter}
\title{Autoregressive Identification of Kronecker Graphical Models}

% \thanksref{footnoteinfo} \thanks[footnoteinfo]{ .}

\author[Padova]{Mattia Zorzi}
\ead{zorzimat@dei.unipd.it}

\address[Padova]{Dipartimento di Ingegneria dell'Informazione, Universit\`a degli studi di
Padova, via Gradenigo 6/B, 35131 Padova, Italy}


\begin{abstract} 
We address the problem to estimate a Kronecker graphical model corresponding to an autoregressive Gaussian stochastic process. The latter is completely described by the power spectral density function whose inverse  has support which admits a Kronecker product decomposition. We propose a Bayesian approach to estimate such a model. We test the effectiveness of the proposed method by some numerical experiments. We also apply the procedure to urban pollution monitoring data.
\end{abstract}

\end{frontmatter}

\section{Introduction}

Graphical models represent a useful tool to describe the conditional dependence structure between Gaussian random variables, \cite{LAURITZEN_1996}. In the present paper we focus on a particular class of graphical models called Kronecker graphical models (KGM), \cite{leskovec2007scalable}. These models received great attention because the corresponding graphs enjoy some properties that emerge in many real graphs, e.g. small diameter and heavy-tailed degree distribution, see \cite{leskovec2010kronecker}. For instance, KGM have been used in recommendation systems 	\citep{allen2010transposable}. {\mz Moreover, KGM can be used to learn basic structures (i.e. modules or groups) useful to  understand the organization of complex networks \citep{leskovec2009networks}.}

In many applications the topology of the graph is not known and has to be estimated from the observed data. \cite{tsiligkaridis2013convergence} consider a KGM corresponding to a Gaussian random vector whose covariance matrix admits a Kronecker product decomposition. Since the graph topology is given by the support of the inverse covariance matrix, the authors proposed a LASSO method for estimating such a graph. However, the assumption that the covariance matrix can be decomposed as a Kronecker product is restrictive in some applications, e.g. this is evident in spatio-temporal MEG/EEG modelling \citep{bijma2005spatiotemporal}. \cite{tsiligkaridis2013covariance} overcame this restriction by considering a Gaussian random vector  whose covariance matrix is a sum of Kronecker products. Moreover, a dynamic extension has been proposed in \cite{8375680}. However, the resulting graphical model is fully connected. Finally,  \cite{CDC_KRON}
considers a KGM corresponding to a Gaussian random vector whose inverse covariance matrix has support which can be decomposed as a Kronecker product. The latter model is less restrictive than the one in  \cite{tsiligkaridis2013convergence} because the covariance matrix does not necessarily admit a Kronecker product decomposition.

The observed signals are typically collected over time and can thus modeled as a high dimensional Gaussian stochastic process. A large body of literature regards the identification of  
sparse graphical models (SGM) corresponding to Gaussian stochastic processes, see  \cite{ARMA_GRAPH_AVVENTI,SONGSIRI_GRAPH_MODEL_2010,SONGSIRI_TOP_SEL_2010,MAANAN2017122,REC_SPARSE_GM,TAC19,alpago2018scalable,zorzi2019graphical}. Such processes are  completely described by the power spectral density (PSD) function. More precisely, the support of the inverse PSD reflects the conditional dependence relations among the components of the process, i.e. the topology of the graph. In all the aforementioned papers the idea is to build a regularized 
maximum likelihood (ML) estimator whose penalty term induces sparsity on the inverse of the PSD. An extension of these models is the introduction of hidden components, see \cite{e20010076,LATENTG,CDC_BRAIN15}. %THREE: \cite{A_NEW_APPROACH_BYRNES_2000,BETA}\\ 
However, the majority of the inference methods for KGM consider i.i.d. processes (i.e. there is no dynamic).

The present paper considers the problem to estimate a KGM 
corresponding to an autoregressive (AR) Gaussian stochastic process. More precisely, we propose a ML estimator adopting a Bayesian perspective. The prior induces the support of the inverse PSD to admit a Kronecker product decomposition. Thus, we do not impose that the PSD admits a Kronecker product decomposition so that the corresponding models is not so restrictive.  {\mz Indeed, if the PSD admits a Kronecker product decomposition, then the dynamic among the nodes in a module is the same in any other module. On the contrary, in our model the dynamic among the nodes in a module is not necessarily the same of those for the other modules.}

In particular, we propose two priors for the ML estimator: the max prior and the multiplicative prior. The latter has been inspired by the ones used in collaborative filtering \citep{yu2009large}, multi-task learning \citep{bonilla2008multi} and it represents the natural extension of the prior proposed in \cite{CDC_KRON} for the the static case. Finally, the penalty term depends on some hyperparameters that we estimate from the data using an {\mz approximate version of the} empirical Bayes approach in the same spirit of \cite{REWEIGHTED}.

The outline of the paper is as follows. In Section \ref{sec:pn_form} we introduce the problem as well as some motivating examples. In Section \ref{sec:additive} we propose the ML estimator for  KGM using the max prior. In Section \ref{sec:additive2} we propose an alternative prior, i.e. the multiplicative prior, to estimate a KGM. In Section \ref{sec:ME} we show that the proposed approach is also connected to a maximum entropy problem. In Section \ref{sec:sim}: we test the proposed methods using synthetic data; we use the method equipped with the max prior to learn the dynamic spatio-temporal graphical model describing the concentration of three urban atmospheric pollutants at a certain area. Finally, in Section \ref{sec:concl} we draw the conclusions.

 \subsection*{Notation}
Given a symmetric matrix $X$, $|X|$ denotes its  determinant, while $X\succ 0$ ($X\succeq 0$) means that $X$ is positive (semi)definite. $X\otimes Y$ denotes the Kronecker product between matrices $X$ and $Y$. Functions on the unit circle $\{e^{i\vartheta} \hbox{ s.t. } \vartheta \in[-\pi,\pi]\}$ will be denoted by capital Greek letters, e.g. $\Phi(e^{i\vartheta})$ with $\vartheta\in[-\pi,\pi]$, and the dependence upon $\vartheta$ will be dropped if not needed, e.g. $\Phi$ instead of $\Phi(e^{i\vartheta})$. If $\Phi(e^{i\vartheta})$ is positive definite (semi-definite) for each $\vartheta\in[-\pi,\pi]$, we will write $\Phi\succ 0$ ($\Phi \succeq 0$). {\mz We denote as $\mathrm{supp}(\Phi)$  the support function of $\Phi$, i.e. the entries of $\Phi$ different from the null function correspond to entries equal to one in $\mathrm{supp}(\Phi)$ otherwise the latter are equal to zero.} The symbol $\Es$ denotes the expectation operator. 
 Given a stochastic process $y=\{\,y(t) ,\; t\in \Zs\,\}$, with some abuse of notation, $y(t)$ will both denote a random vector and its sample value. {\mz The notation $A\,\bot\, B\,|\,X$ means that the vector subspaces $A$ and $B$ of a Hilbert space are conditionally orthogonal given a third subspace $X$.}

\section{Problem formulation}\label{sec:pn_form}
\begin{figure}[htbp]
\centering
  \includegraphics[width=0.7\columnwidth]{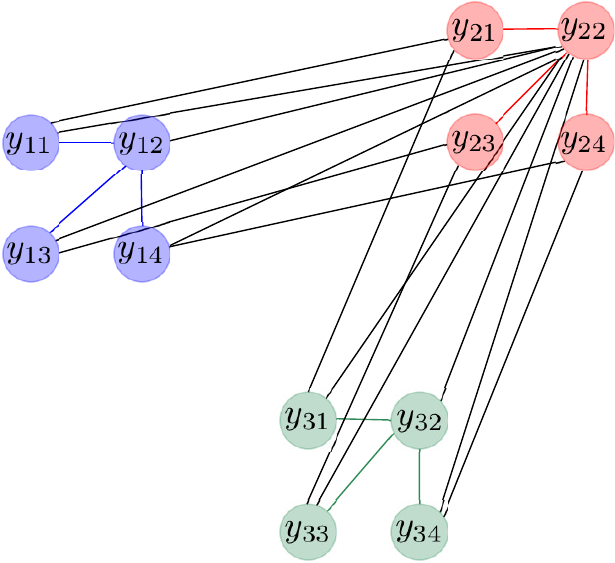}  
\caption{An example of dynamic KGM with $m_1=3$ modules (red, blue and green) composed by $m_2=4$ nodes. $y=[\, y_{11}\, y_{12}\, y_{13}\, y_{14}\, y_{21}\, y_{22}\, y_{23}\, y_{24}\, y_{31}\, y_{32}\, y_{33}\, y_{34}\,]^T$, $\Ec_1=\{(1,2),(2,3)\}$ and $\Ec_2=\{(1,2),(2,3),(2,4)\}$.} \label{fig:graph_ex}
\end{figure}

Consider an AR Gaussian discrete-time zero mean full rank stationary stochastic process denoted by  $y=\{\, y(t),\;  t\in \mathbb{Z}\}$ where $y(t)\in\Rs^{m_1 m_2}$, $m_1,m_2\in\Ns$. Such a process is completely characterized by its PSD
\al{\Phi(e^{i\vartheta})=\sum_{s\in \mathbb{Z}}  e^{-is \vartheta} R_s, \; \; \; \vartheta\in[-\pi,\pi]} where $R_s=\Es[ y(t)y(t+s)^T]$, with $s\in \mathbb{Z}$.  Notice that $\Phi^{-1}\in \Qc_{m_1m_2,n}$ where \al{ \label{set_Qc_mn_reparametrized}\Qc_{m_1m_2,n}=\left\{\, \sum_{s=-n}^n Q_se^{-i\theta s} \hbox{ s.t. } Q_s=Q_{-s}^T\in \Rs^{m_1 m_2}\right\}} is the family of pseudo-polynomial matrices and $n\in\Ns$ denotes the order of the AR process. Such a model admits the following interpretation in terms of a dynamic graphical model describing conditional dependence relations \citep{REMARKS_BRILLINGER_1996}. Let $y_{hk}$ denote the entry of $y$ in position $(h-1)m_2+k$ with $h\in \Vc_1$ and $k\in\Vc_2$ where $\Vc_1=\{ 1\ldots m_1  \}$ and $\Vc_2=\{ 1\ldots m_2  \}$. Given $I\subseteq \Vc_1\times \Vc_2$, we denote as 
 \al{\chi_{I}=\overline{\mathrm{span}}\{\, y_{hk}(t) \hbox{ s.t. }\, (h,k)\in I\, , \, t\in \Zs\, \}}
the closure of all finite linear combinations of $y_{hk}(t)$ with $h\in I_1$, $k\in I_2$ and $t\in \Zs$. {\mz The latter is a vector subspace of the Hilbert space of Gaussian random variables having finite second order moments. Let $(h,k)\neq (j,l)$, then $y_{hk}$ and $y_{jl}$ are conditionally independent if and only if \al{\chi_{\{ (h, k)\}} \, \bot\, \chi_{\{(j, l)\}} \, |\, \chi_{\Vc_1 \times \Vc_2  \setminus \{(h,k),(j,l)\}},}
see Section 2 in \cite{LINDQUIST_PICCI} for more details.} We assume these conditional dependence relations define a dynamic KGM $\Gc(\Vc_1\times \Vc_2 ,\Ec_1\times \Ec_2)$ where $\Vc_1\times \Vc_2$ and {\mz $\Ec_1\times \Ec_2$} denote the set of nodes and edges, respectively, with $\Ec_1\subseteq \Vc_1 \times \Vc_1$ and $\Ec_2\subseteq \Vc_2 \times \Vc_2$. More precisely, the nodes represent the components $y_{hk}$ of $y$ and the lack of an edge in $\Ec_1$ or in $\Ec_2$ means conditional independence:
\al{ \label{cond_ind_cond}&(h,j) \notin \Ec_1 \hbox{ or } (k,l)\notin \Ec_2\, \, \iff\nn \\ 
& \hspace{2cm} \chi_{\{ (h, k)\}} \, \bot\, \chi_{\{(j, l)\}} \, |\, \chi_{\Vc_1\times \Vc_2  \setminus \{(h,k),(j,l)\}}.}
In graph $\Gc$ we can recognize $m_1$ modules containing $m_2$ nodes and sharing the same graphical structure described by $\Ec_2$, while the interaction among those $m_1$modules is described by $\Ec_1$. An example of dynamic KGM is provided in Figure \ref{fig:graph_ex}.
 \cite{Dahlhaus2000} showed that conditional dependence relations are characterized by the support of $\Phi^{-1}$. Therefore, in our case an equivalent condition of (\ref{cond_ind_cond})  is 
 \al{\label{cond_sup} \mathrm{supp}(\Phi^{-1})= E_1\otimes E_2}
 where $E_1$ and $E_2$ are {\mz adjacency} matrices of dimension $m_1\times m_1$ and $m_2 \times m_2$, respectively, such that $(E_1)_{hj}=1$ if and only if $(h,j)\in\Ec_1$ otherwise 
 $(E_1)_{hj}=0$ and likewise for $E_2$ and $\Ec_2$. In other words, $y$ corresponds to a dynamic KGM if and only if condition (\ref{cond_sup}) holds. The next proposition shows that $\Ec_1$ and $\Ec_2$ describe the conditional dependence relations among modules and nodes in each module, respectively. 
  \pro \label{prop_MEG} Consider the vector spaces 
 \al{\chi^\star_{I_1}&=\overline{\mathrm{span}}\{\, y_{hk}(t) \hbox{ s.t. }\, h\in I_1, \,k\in \Vc_2   , \, t\in \Zs\, \}\\
 \chi^\dag_{I_2}&=\overline{\mathrm{span}}\{\, y_{hk}(t) \hbox{ s.t. }\, h\in \Vc_1,\, k\in I_2   , \, t\in \Zs\, \} } where $I_1\subseteq \Vc_1$ and $I_2\subseteq \Vc_2$. Under assumption (\ref{cond_sup}) we have that 
 \al{& \label{cond_MEG1}(h,j)\notin \Ec_1\; \iff\; \chi_{ \{h\}}^\star \, \bot\, \chi_{\{ j\}}^\star \, |\, \chi^\star_{\Vc_1 \setminus \{h,j\} }\\
& \label{cond_MEG2} (k,l)\notin \Ec_2\; \iff\;\chi_{ \{ k\}}^\dag \, \bot\, \chi_{\{ l\}}^\dag \, |\, \chi^\dag_{\Vc_2 \setminus \{k,l\} }.}
 \epro
 \proof We prove condition (\ref{cond_MEG1}); condition (\ref{cond_MEG2}) can be proved in a similar way. Let 
 \al{z_h(t)=[\, y_{h1}(t) \ldots y_{hm2}(t)\,]^T.} 
 Let $\varepsilon(t)=[\,\varepsilon_h(t)^T\, \varepsilon_j(t)^T\,]^T$ be the projection error of $[\,z_h(t)^T\, z_j(t)^T\,]^T$ onto $\chi^\star_{\Vc_1\setminus \{h,j\}}$ with $h\neq j$. Since $\varepsilon$ is a zero mean Gaussian stationary process, proving (\ref{cond_MEG1}) is equivalent to prove that 
 \al{\label{exp_errr}(h,j)\notin \mathcal E_1\, \iff\,\Es[\varepsilon_h(t)\varepsilon_j(t)^T]=0,\; \; \forall\, t\in\Zs.}
 Let $P_\Pi$ be the permutation matrix that permutes the components of $y$ in order to obtain
 \al{\tilde y=P_\Pi y= \left[\begin{array}{c}z_h \\ z_j\\ \hline z_{k\neq h,i} \end{array}\right]=\left[\begin{array}{c}\tilde y_1 \\ \hline \tilde y_2\\  \end{array}\right]}  where $z_{k\neq h,j} $ is the process obtained by stacking $z_k$ with $k\neq h,j$. We partition the PSD $\tilde \Phi$ of $\tilde y$ in conformable way 
 \al{\tilde \Phi= \left[\begin{array}{cc}  \tilde \Phi_{1,1} & \tilde \Phi_{1,2} \\ \tilde \Phi_{2,1}  & \tilde \Phi_{2,2} \\  \end{array}\right].} Notice that \al{\tilde \Phi^{-1}= P_\Pi \Phi^{-1}P_\Pi^T=\left[\begin{array}{cc}  \Sigma & \hspace{0.3cm}\star\hspace{0.2cm} \\ \star  &\hspace{0.3cm} \star\hspace{0.2cm} \\  \end{array}\right]}
 where \al{\Sigma^{-1}=  \tilde \Phi_{1,1} - \tilde \Phi_{1,2}  \tilde \Phi_{2,2}^{-1} \tilde \Phi_{2,1}=\tilde \Phi_{\varepsilon} }
and  $\tilde \Phi_\varepsilon$ denotes the PSD of $\varepsilon$. Therefore, $\Es[\varepsilon_h(t)\varepsilon_j(t)^T]=0$ $\forall\, t\in\Zs$ if and only if $\Sigma$ is block-diagonal (according to the partition in $\varepsilon(t)=[\,\varepsilon_h(t)^T\, \varepsilon_j(t)^T\,]^T$). In view of (\ref{cond_sup}), $\Sigma$ is block-diagonal  if and only if $(h,j)\notin \Ec_1$. We conclude that  (\ref{exp_errr}) holds. \qed\\

 Throughout the paper we want to address the following identification problem about dynamic KGM. 
 \pb \label{ref_pb}Consider an AR Gaussian zero mean full rank stationary stochastic process $y$ of order $n$ and taking values in $\Rs^{m_1 m_2}$. Assume that $m_1$, $m_2$, $n$ are known and collect a finite length sequence $ y^N:=\{y(1), y(2) \ldots y (N)\}$ extracted from a realization of $y$. Let $\Phi$ be the PSD of $y$ satisfying (\ref{cond_sup}). Find an estimate $\hat \Phi$ of $\Phi$ from $y^N$ such that $\mathrm{supp}(\hat \Phi^{-1})= \hat E_1\otimes \hat E_2$ where $\hat E_1$ and $\hat E_2$ represent an estimate of $E_1$ and $E_2$, respectively.
 
 \epb It is worth noting that condition (\ref{cond_sup}) is weaker than $\Phi^{-1}$ admit a Kronecker decomposition, i.e. $\Phi^{-1}=\Phi^{-1}_1 \otimes \Phi^{-1}_2$. {\mz Accordingly, we do not constrain the dynamics in each module to be same.} In what follows we present some practical problems in which a stochastic process corresponding to a dynamic KGM could be used. 
 
\subsection{Dynamic spatio-temporal modeling}
Consider a non-stationary zero mean Gaussian process {\mz $x=\{x(\tilde t), \; \tilde t\in\Zs\}$ indexed by the time variable $\tilde t$ where $x(\tilde  t)$} takes values in $\Rs^{m_2}$, $m_2\in\Ns$, and whose covariance lags sequence {\mz $P_{\tilde t,\tilde s}=\Es[x(\tilde t) x(\tilde t+\tilde s)^T]$ is such that 
\al{P_{\tilde t+m_1,\tilde s}=P_{\tilde t,\tilde s},\;\;\; \tilde t,\tilde s\in\Zs}
with $m_1\in\Ns$.} We can rewrite $x$ in terms of the stochastic process $y=\{y(t), \; t\in\Zs\}$ defined as
\al{y(t)=[\, x((t-1)m_1+1)^T \ldots x(tm_1)^T\,]^T.} {\mz It is worth noting that the time variables $\tilde t$ and $t$ are different: there is a decimation relationship between them and the decimation factor is $m_1$.} It is not difficult to see that $y$ is zero mean, Gaussian and stationary. In particular, its covariance lags sequence $R_s=\Es[y(t)y(t+s)^T]$ is {\tiny \al{R_s=\left[\begin{array}{cccccc}P_{1,sm_1} & P_{1,sm_1+1} & \ldots & \ldots &\ldots &  P_{1,(s+1)m_1-1} \\ P_{1,sm_1+1}^T & P_{2,sm_1}^T & \ddots & &  & \vdots\\
\vdots & \ddots & \ddots & & & \vdots \\ \vdots &   &  &   & &  P_{m_1-1,sm_1+1}\\ P_{1,(s+1)m_1-1}^T & \ldots  & \ldots  & \ldots & P_{m_1-1,sm_1+1}^T  & P_{m_1,sm_1}
\end{array}\right]\nn} }
where we have exploited the relation {\mz $P_{\tilde t,\tilde s}^T=P_{\tilde t+\tilde s,-\tilde s}$}.
Let $\Phi$ be the PSD of $y$ and assume that condition (\ref{cond_sup}) holds. In view of Proposition \ref{prop_MEG}, with $y_{hk}(t)=x_k((t-1)m_1+h)$, $\Ec_2$ describes the conditional dependence relations among the components $x_k$, with $k\in\Vc_2$, of process $x$.  We conclude that $x$ can be understood as a spatio-temporal process.   In the special case that 
 {\mz \al{\label{cond_MEG}P_{\tilde t,\tilde s}=0,\; \; \forall \, \tilde s \hbox{ s.t. }|\tilde s|\geq m_1}} we have that $R_s=0$ for any $s>0$ and thus $y$ is an i.i.d. Gaussian process, i.e. $\Phi$ is a covariance matrix. The latter  models magnetoencephalography (MEG) measurements used for mapping brain activity, \cite{bijma2005spatiotemporal}. More precisely, $z_k(t)=[\, x_k((t-1)m_1+1)\,\ldots\, x_k(tm_1) \,]^T$ models the measurements at the $k$-th brain area during the $t$-th trial of length $m_1$. All the trials are independent. Moreover, the latter are identically distributed because in each trial the patient is required to perform the same cognitive
task. In our framework we can remove condition (\ref{cond_MEG}), i.e. the trials now can be dependent. This means 
that in the future such trials can be scheduled in a sequential way and modeled through the spatio-temporal process $x$. In view of Proposition \ref{prop_MEG}, $\Ec_2$ describes the conditional dependence relations between the different brain areas.
  It is worth noting that we do not force the Kronecker structure on $\Phi^{-1}$: such freedom has shown to be crucial for an effective MEG modeling. Finally, the spatio-temporal process $x$ can be potentially used also to urban pollution monitoring, see Section \ref{sec:sec:poll}.

\subsection{Multi-task modeling}

We consider a network composed by $m_2\in\Ns$ agents (i.e. nodes); each agent is described by a zero mean Gaussian stationary stochastic process.  We want to model such a network under $m_1\in\Ns$ heterogeneous conditions (i.e. tasks). Let $y_{hk}=\{y_{hk}(t), \; t\in \Zs\}$ denote the stochastic process describing the $k$-th agent under the $h$-th task. Then, we can model the network through the $m_1$ stochastic processes 
\al{ \label{all_prc_MT}w_h=[\, y_{h 1}\,\ldots\, y_{h m_2}  \,]^T,\; \; \; h\in\Vc_1}      
where $h$ denotes the task. A more flexible approach is to model all processes in (\ref{all_prc_MT}) together in order to exploit commonalities and differences across the tasks, see \cite{allen2010transposable,yu2009large}. More precisely, we consider the stationary stochastic process $y$ obtained by stacking $w_h$ with $h\in\Vc_1$. Let $\Phi$ be {\mz the} PSD  of $y$ 
and such that (\ref{cond_sup}) holds. In view of Proposition \ref{prop_MEG}, $\Ec_1$ describes the conditional dependence relations among the tasks, while $\Ec_2$ describes the ones among the agents of the network. It is worth noting that our model is dynamic in contrast with the ones in  \cite{allen2010transposable,yu2009large} which are static. The proposed model could be used to describe the travel demand in a public transport system of a certain city network, \cite{chidlovskii2017multi}. More precisely, $y_{hk}(t)$ denotes the total number of boarding events at day $t$, $k$ represents a particular area in the city network and $h$ represents the type of transportation (e.g. bus, train, tram).

\section{Identification of KGM} \label{sec:additive}
We aim to solve Problem \ref{ref_pb} where we parametrize the PSD of $y$ as $\Phi=\Sigma^{-1}$ where $\Sigma\in \Qc_{m,n}$, $m=m_1m_2$ and 
{\mz \al{ \Sigma(e^{i\theta}) =S_0+\frac{1}{2}\sum_{t=1}^n [ S_t e^{-it\vartheta}+S_t^T e^{it\vartheta}].}} Notice that 
{\mz  \al{&(\Sigma(e^{i\vartheta}))_{hk,jl}\nn \\ &=(S_0)_{hk,jl}+\frac{1}{2}\sum_{t=1}^n [(S_t)_{hk,jl}e^{-i t \vartheta}+(S_t)_{jl,hk}e^{i t \vartheta}]\nn}}
where $(\Sigma)_{hk,jl}$ denotes the entry of $\Sigma$ with row $(h-1)m_2+k$ and column $(j-1)m_2 +l$ with $h,j\in\Vc_1$ and $k,l\in\Vc_2$; the same meaning has the notation $(S_t)_{hk,jl}$ for $S_t$. We consider the following regularized ML estimator of $\Sigma$, and thus of $\Phi$: 
\al{\label{reg_ml}\hat \Sigma= &\underset{\Sigma\in \Qc_{m,n}}{ \mathrm{argmin}} \,\ell( y^N; \Sigma)+  g(\Sigma;\Lambda,\Gamma)\nn\\
& \hbox{ s.t. }  \Sigma\in\Qc_{m,n}^+}
where $\Qc_{m,n}^+:=\{\Sigma\in\Qc_{m,n} \hbox{ s.t. } \Sigma\succ 0\}$.
The term $\ell( y^N; \Sigma)$ is an approximation of the negative log-likelihood of $  y(n
+1) \ldots  y(N)$ given $  y(1) \ldots  y (n)$ under the assumption that $y$ is an AR process of order $n$: 
\al{\label{cond_PDF}\ell(  y^N; \Sigma)= \frac{N-n}{4\pi}\int_{-\pi}^\pi [-\log |\Sigma| + \tr(\hat \Phi_{p}\Sigma)]\mathrm d\vartheta+c} 
where  \al{ \label{windowed_cor}  \mz \hat  \Phi_{p}&=\mz \hat R_0+\frac{1}{2}\sum_{s=1}^n  [\hat R_s e^{-is\vartheta}+\hat R_s^T e^{is\vartheta}]\nn \\
\hat R_s & =\frac{1}{N-n} \sum_{t=1}^{N-s} y(t)y(t+s)^T
} and $c$ is a term not depending on $y^N$ and $\Sigma$. Notice that $\hat R_s$ represents an estimate of $R_s$ from data $ y^N$ and $\hat \Phi_{p}$ is the truncated periodogram of $\Phi$. Let $\Td(\hat R)$ denote the block-Toeplitz matrix whose first block row is $[\,\hat R_0 \, \hat R_1 \ldots \hat R_n\,]$. Throughout the paper we make the assumption that $\Td(\hat R)\succ 0$. The latter assumption holds for $N$ sufficiently large since $y$ is a full rank process. {\mz The penalty  term $g(\Sigma;\Lambda,\Gamma)$ induces some desired properties in the solution $\hat\Sigma$; $\Lambda$ and $\Gamma$ are the regularization matrices (hereafter called hyperparameter matrices) that later will be estimated from the data.}  \cite{REWEIGHTED} proposed the following penalty term{\mz , which in turn is built on the one by \cite{SONGSIRI_TOP_SEL_2010},} for estimating a SGM:
\al{\label{hSPARSE}g (\Sigma;&\Omega) =\sum_{(h,k,j,l)\in\Tc_S}  (\Omega)_{hk,jl} \times \nn \\ &\times \max\{\max_{t=0\ldots n} |(S_t)_{hk,jl}|,\max_{t=0\ldots n} |(S_t)_{jl,hk}| \}}
where $\Tc_S=\{(h,k,j,l) \hbox{ s .t. } (h-1)m_2+k\geq (j-1)m_2+l\}$ and the entries of the weight symmetric matrix $\Omega$ are estimated using an {\mz approximate version of the} empirical Bayes approach. It has been shown that (\ref{hSPARSE}) leads to an estimate of $\Phi$ whose inverse is sparse. Here, instead, we consider a penalty term  $g(\,\cdot\,; \Lambda,\Gamma)\, :\, \Qc_{m,n} \longrightarrow \Rs$ which  is designed in such a way to induce (\ref{cond_sup}) in the solution of (\ref{reg_ml}). More precisely, we consider  
 \al{\label{hA}   g &(\Sigma;\Lambda,\Gamma) =\hspace{-0.2cm}\sum_{(h,k,j,l)\in\Tc}  \hspace{-0.4cm} \max\{\lambda_{hj},\gamma_{kl}\}q_{hk,jl}(\Sigma)}  
with\al{ & \Tc :=\{(h,k,j,l) \hbox{ s.t. } h,j\in\Vc_1 \; k,l\in\Vc_2 \hbox{ and } h\geq j,\;  k\geq l  \}\nn\\
 & q_{hk,jl}(\Sigma) := \max\{  \max_{t=0 \ldots n} |(S_t)_{hk,jl}|, \max_{t=0 \ldots n} |(S_t)_{hl,jk}|,\nn\\ 
 &\hspace{2.7cm}\max_{t=0 \ldots n} |(S_t)_{jl,hk}|, \max_{t=0 \ldots n} |(S_t)_{jk,hl}|\}, \nn}   
 $\lambda_{hj},\gamma_{kl}\geq 0$, with $h,j\in\Vc_1$ and $k,l\in\Vc_2$. \begin{figure}[htbp]
\centering
  \includegraphics[width=0.7\columnwidth]{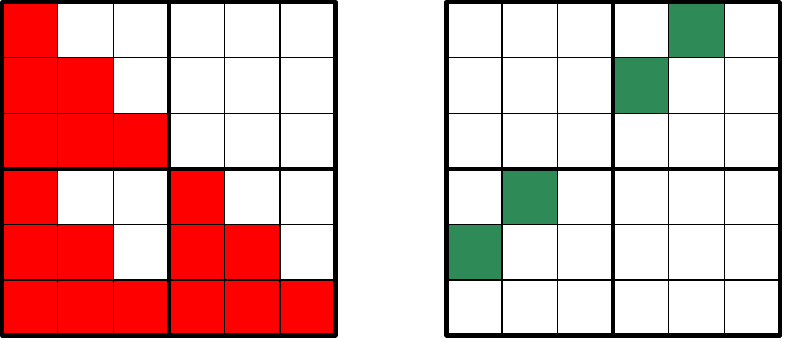}  
\caption{\mz Example with $m_1=2$ and $m_2=3$. {\em Left.} In red the positions of $S_t$ which are taken in $\mathcal T$. {\em Right}.  In green the entries of $S_t$ which are considered in $q_{hk,jl}(\Sigma)$ with $h=2$, $j=1$, $k=2$, $l=1$.}\label{fig:set2}
\end{figure}  {\mz It is worth noting that the index set $\mathcal T$ contains only a subset of all the possible indexes characterizing the entries of $\Sigma$, see Figure \ref{fig:set2} ({\em left}). This is because the support of $\Sigma$ has to satisfy the symmetric Kronecker  structure in (\ref{cond_sup}), indeed recall that $E_1$ and $E_2$ are symmetric adjacency matrices. Accordingly, we  induce a group sparsity not only to guarantee that $St$'s have the same support but also to guarantee the symmetric property in (\ref{cond_sup}), see the example of Figure \ref{fig:set2} (right).} $\lambda_{hj}$ and $\gamma_{kl}$ are the entries in position $(h,j)$ and $(k,l)$ of the symmetric matrices $\Lambda$ and $\Gamma$, respectively. Notice that $\Lambda$ and $\Gamma$ are nonnegative matrices of dimension $m_1\times m_1$ and $m_2\times m_2$, respectively. Some comments about this penalty term follow. $g(\Sigma;\Lambda,\Gamma)$ is a weighted sum of the nonnegative terms $q_{hk,jl}(\Sigma)$, so the penalty  induces many of these terms to be equal to zero. If $q_{hk,jl}(\Sigma)=0$, then $(\Sigma)_{hk,jl}$, $(\Sigma)_{jl,hk}$, $(\Sigma)_{hl,jk}$ and $(\Sigma)_{jk,hl}$ coincide with the null function. Therefore, such a penalty encourages a common sparsity pattern {\mz (i.e. group sparsity, see as example \cite{DUE})} on $\Sigma$ according to (\ref{cond_sup}). More precisely, if $\lambda_{hj}$ is large then it is more likely that the solution of (\ref{reg_ml}) is such that $q_{hk,jl}(\Sigma)=0$ for any $k,l\in\Vc_2 $, i.e. (\ref{cond_sup}) holds with $(h,j)\notin \Ec_1$. If $\gamma_{kl}$ is large then it is more likely that the solution of (\ref{reg_ml}) is such that $q_{hk,jl}(\Sigma)=0$ for any $h,j\in\Vc_1 $, i.e. (\ref{cond_sup}) holds with $(k,l)\notin \Ec_2$.

 \pro\label{prop_reg_MLA} The objective function of Problem (\ref{reg_ml}) is strictly convex in $\Sigma$. Moreover, Problem (\ref{reg_ml}) {\mz admits} a unique solution.  
 \epro
 \proof  It is not difficult to prove that $\ell(y^N; \Sigma)$ is strictly convex over $\Qc_{m,n}^+$, see for instance Theorem 5.1 in \cite{FERRANTE_TIME_AND_SPECTRAL_2012}. Moreover $g(\Sigma;\Lambda,\Gamma)$ is a {\mz convex} function of $\Sigma$, accordingly the objective function in (\ref{reg_ml}) is strictly convex. Then, it is possible to prove that (\ref{reg_ml}) {\mz admits} solution (which is clearly in the interior of $\Qc_{m,n}^+$) using arguments similar to the ones of Theorem 1 in \cite{REWEIGHTED}. Finally, the uniqueness of such a solution follows from the strict convexity of the objective function. \qed\\
 As we already noticed, the support of $\hat \Sigma$ depends on the weights $\lambda_{hj}$ and $\gamma_{kl}$. {\mz The latter can be estimated from the data.} With this aim we adopt the following Bayesian perspective: the parameters $S_t$, $t=0\ldots n$ characterizing $\Sigma$ are random matrices with a suitable probability density function (PDF) or simply prior, i.e. $\Sigma$ is a random process taking values in $\Qc_{m,n}^+$. We denote by $p(\Sigma|\Lambda,\Gamma)$ the PDF of $\Sigma$ given $\Lambda$ and $\Gamma$. We assume that $(\Sigma)_{hk,jl}$'s are independent {\mz of} each other, so that 
\al{\label{form_indip}p(\Sigma|\Lambda,\Gamma)=\prod_{(h,k,j,l)\in\Tc} p ((\Sigma)_{hk,jl}|\lambda_{hj},\gamma_{kl}). } 
 Moreover, we define the ``max prior'' \al{ \label{prior_SIGMA}p((\Sigma)_{hk,jl}|\lambda_{hj},\gamma_{kl})=\frac{e^{-\max\{\lambda_{hj},\gamma_{kl}\}q_{hk,jl}(\Sigma)}}{c_{hk,jl}} }
where $c_{hk,jl}>0$ is the normalizing constant. {\mz It is worth noting that (\ref{prior_SIGMA}) is a Gibbs prior; such a choice is very common in Bayesian methods because it leads to a solvable problem, see e.g. \cite{UNO}.} Therefore, we have
\al{\label{p_SigmaA} p(\Sigma|\Lambda,\Gamma)=\frac{e^{- g (\Sigma;\Lambda,\Gamma)}}{\prod_{(h,k,j,l)\in\Tc} c_{hk,jl}}.} $\Lambda$ and $\Gamma$ are modeled as independent random matrices whose entries are independent {\mz of} each other with exponential distribution:
\al{&  \label{form_LG_1}p(\Lambda,\Gamma)=p(\Lambda)p(\Gamma) \\
& \label{form_LG_2} p(\Lambda)=\prod_{h\geq j} p(\lambda_{hj}), \; \;\;  p(\lambda_{hj})= \varepsilon e^{-\varepsilon \lambda_{hj}} \\ 
& \label{form_LG_3}p(\Gamma)=\prod_{k\geq l} p(\gamma_{kl}), \; \;\;  p(\gamma_{kl})= \varepsilon e^{-\varepsilon \gamma_{kl}}  }
where $\varepsilon>0$ is a fixed small constant.

{\mz \rema We assume the elements of $\Lambda$ are independent because there is no a priori information on how the modules are related  (e.g. an a priori information like: ``if there is an edge between the modules $h$ and $j$, then it is likely that there is a link between the modules $\tilde h$ and $\tilde j$''). Likewise, the elements of $\Lambda$ are independent because there is no a priori information on how the nodes in each module are related. Finally, the elements of $\Lambda$ and $\Gamma$ are independent because there is no a priori information on how a module and a node in a module are related.\erema \rema It is worth noting that the Bayesian model introduced in this section is only used to estimate the KGM, indeed we do not assume that the actual model is generated in this way. The same perspective has been used in \cite{chen2012estimation}.\erema}
  
The negative log-likelihood of $y^N$ and $\Sigma$ takes the form:
\al{\ell ( y^N,&\Sigma,\Lambda,\Gamma)=-\log p( y^N,\Sigma,\Lambda,\Gamma)\nn\\
& = -\log p( y^N|\Sigma)-\log p(\Sigma|\Lambda,\Gamma) -\log p(\Lambda)p(\Gamma)\nn.}
The negative log-conditional PDF $-\log p( y^N|\Sigma)$  {\mz coincides} with (\ref{cond_PDF}), thus 
\al{\label{ML_completa}\ell ( y^N,\Sigma,\Lambda,\Gamma)=\ell( &y^N; \Sigma)+g(\Sigma;\Lambda,\Gamma)+\hspace{-0.4cm}\sum_{(h,k,j,l)\in\Tc}\hspace{-0.4cm} \log c_{hk,jl} \nn\\
&+ \sum_{h\geq j} \varepsilon \lambda_{hj}+ \sum_{k\geq l} \varepsilon \gamma_{kl}} where we discarded the terms not depending on $\Sigma$, $\Lambda$ and $\Gamma$.

\pro \label{pro_C_somma} Consider the prior of $\Sigma$ given by (\ref{prior_SIGMA}). Then, we have \al{ c_{hk,jl}\leq \left\{\begin{array}{ll} 
\upsilon_{hk,jl}(\lambda_{hj}\gamma_{kl})^{-(n+1)}, & \hbox{ $h=j$ and $k=l$} \\ 
\upsilon_{hk,jl}(\lambda_{hj}\gamma_{kl})^{-(2n+1)}, & \hbox{  $h=j$ and $k>l$} \\ 
\upsilon_{hk,jl}(\lambda_{hj}\gamma_{kl})^{-(2n+1)}, & \hbox{ $h>j$ and $k=l$} \\ 
\upsilon_{hk,jl}(\lambda_{hj}\gamma_{kl})^{-(4n+2)}, & \hbox{  $h>j$ and $k>l$}  
\end{array}\right.} 
where $\upsilon_{hk,jl}$ are constant terms not depending on $\Lambda$ and $\Gamma$. \epro

\proof Since $c_{hk,jl}$ is the normalizing constant in (\ref{prior_SIGMA}), we have 
\al{c_{hk,jl}&=\int_{\Qc_{m,n}^+} e^{-\max\{\lambda_{hj},\gamma_{kl}\}q_{hk,jl} (\Sigma)}\mathrm d \Sigma\nn\\
& \leq \int_{\Qc_{m,n}} e^{-\max\{\lambda_{hj},\gamma_{kl}\}q_{hk,jl} (\Sigma)}\mathrm d \Sigma \nn\\
&\leq \upsilon_{hk,jl}(\max\{\lambda_{hj},\gamma_{kl}\})^{-\mathrm n \sharp(q_{hk,jl} (\Sigma))}}
where $\mathrm n \sharp(q_{hj,kl} (\Sigma))$ is the number of parameters characterizing $(q_{hk,jl} (\Sigma))$ and $\upsilon_{hj,kl}$ is a term not depending on $\Lambda$, $\Gamma$ and $\Sigma$, see Lemma 9 in \cite{REWEIGHTED}. For instance, in the case that $h=j$ and $k=l$, we have 
\al{q_{hk,hk}(\Sigma)= \max_{t=0\ldots n} |(S_t)_{hk,hk}|} which depends on $n+1$ parameters. 
\qed\\
In view of Proposition \ref{pro_C_somma}, we have an upper bound for $\ell(y^N, \Sigma, \Lambda,\Gamma)$:
\al{\label{upper_bound_somma}\ell( & y^N, \Sigma, \Lambda,\Gamma)\leq \tilde \ell(  y^N, \Sigma, \Lambda,\Gamma) := \ell(y^N; \Sigma)\nn\\ &+g(\Sigma;\Lambda,\Gamma)-\hspace{-0.4cm} \sum_{(h,k,j,l)\in\Tc} \hspace{-0.4cm} \alpha_{hk,jl} \log \max\{\lambda_{hj},\gamma_{kl}\}\nn\\ &+\varepsilon\sum_{h\geq j} \lambda_{hj}+\varepsilon\sum_{k\geq l} \gamma_{kl}}
where \al{\alpha_{hk,jl}= \left\{\begin{array}{ll}
n+1 , & \hbox{if $h=j$ and $k=l$} \\ 
 2n+1 , & \hbox{if $h=j$ and $k>l$} \\ 
 2n+1 , & \hbox{if $h>j$ and $k=l$} \\ 
 4n+2 , & \hbox{if $h>j$ and $k>l$}.\end{array}\right.\nn} 
 Following the generalized maximum likelihood (GML) method \citep{zhou1997approximate}, an estimator for $\Sigma$, $\Lambda$ and $\Gamma$ is given by 
\al{\label{pb_GML}(\hat \Sigma,\hat \Lambda,\hat \Gamma) =& \underset{ (\Sigma,  \Lambda,  \Gamma)\in\Cc }{\;\mathrm{argmin}}\; \tilde \ell ( y^N,  \Sigma ;\Lambda,  \Gamma )}
where \al{\Cc:=\{ (\Sigma,  \Lambda,  \Gamma)& \hbox{ s.t. } \Sigma\in\Qc_{m,n}^+, \; \lambda_{hj}\geq 0,\;\gamma_{kl}\geq 0
\}.} 

\rema It is worth noting that an estimator of $\Lambda$ and $\Gamma$ can be computed by maximizing the marginal PDF $p(y^N)$ under model (\ref{p_SigmaA})-(\ref{form_LG_3}), see \cite{friedman2001elements}. On the other hand, it is not possible to find an analytical expression for the latter. 
\erema
\pro \label{prop_ex_GML}Problem (\ref{pb_GML}) {\mz admits} solution.\epro
\proof First, note that $\Cc$ is open and unbounded. We show that   (\ref{pb_GML}) is equivalent to the minimization of $\tilde \ell $ over a compact set $\Cc^\star$. To this aim, we consider the set 
\al{\bar \Cc:=\{(\Sigma,  \Lambda,  \Gamma)& \hbox{ s.t. } \Sigma\in\bar \Qc_{m,n}^+,  \; \lambda_{hj}\geq 0,\;\gamma_{kl}\geq 0\}} where $\bar \Qc_{m,n}^+$ denotes the closure of  $ \Qc_{m,n}^+$. In Lemma 5.1 \citep{FERRANTE_TIME_AND_SPECTRAL_2012} it has been shown that it is possible to extend $\ell(y^N; \Sigma)$ over $\bar Q_{m,n}^+$; more precisely, $\ell(y^N;\Sigma)$ is lower semicontinuous in $\bar Q_{m,n}^+$ with values in the extended reals. Therefore, 
we can extend  $\tilde \ell(y^N, \Sigma,\Lambda,\Gamma)$ over $\bar \Cc$. The latter  is lower semicontinuous in $\bar \Cc$ with values in the extended reals. Next, we show   that sequences $ (\Sigma^{(r)},  \Lambda^{(r)},  \Gamma^{(r)})$ which diverge or approach the boundary of $\bar \Cc$ cannot be infimizing sequences of $\tilde \ell$. More precisely, we may have four cases.  \\
$\bullet$ {\em Case 1:}  $(\Sigma^{(r)},\Lambda^{(r)},\Gamma^{(r)})\in \Cc$ is a convergent sequence such that  there exists at least one $(\bar h,\bar j)$ and/or one $(\bar k,\bar l)$ for which $\lambda_{\bar h,\bar j }^{(r)}\rightarrow \infty$ and/or $\gamma_{\bar k,\bar l }^{(r)}\rightarrow \infty$. Using the fact that $g(\Sigma; \Lambda,\Gamma)$ is nonnegative, we have 
\al{&\lim_{r\rightarrow \infty}  \tilde \ell(y^N,\Sigma^{(r)},\Lambda^{(r)},\Gamma^{(r)})\geq \lim_{r\rightarrow \infty} \ell(y^N,\Sigma^{(r)})\nn\\ &-\hspace{-0.4cm}\sum_{(h,k,j,l)\in\Tc} \hspace{-0.4cm} \alpha_{hk,jl} \log\max\{\lambda_{hj}^{(r)},\gamma_{kl}^{(r)}\} +\varepsilon \sum_{h\geq j} \lambda_{hj}^{(r)}+\varepsilon \sum_{k\geq l} \gamma_{kl}^{(r)}\nn \\ 
&  \geq c +\lim_{r\rightarrow \infty} -\hspace{-0.4cm} \sum_{(h,k,j,l)\in\Tc}\hspace{-0.4cm}  \alpha_{hk,jl} \log\max\{\lambda_{hj}^{(r)},\gamma_{kl}^{(r)}\}\nn\\ &+ \varepsilon \sum_{h\geq j}   \lambda_{ h  j}^{(r)}  + \varepsilon\sum_{k\geq l}  \gamma_{kl}^{(r)}=\infty \nn}  
 where we exploited the following facts: $\ell(y^N;\Sigma)$ is bounded from below on $\bar \Qc_{m,n}^+$, i.e. there exists a finite constant $c$ such that $\ell(y^N;\Sigma)\geq c $ for any $\Sigma\in\bar \Qc_{m,n}^+$, see Lemma 5.3 in \cite{FERRANTE_TIME_AND_SPECTRAL_2012}; the logarithmic term $\alpha_{hk,jl} \log\max\{\lambda_{hj},\gamma_{kl}\}$ is dominated by the linear term $\varepsilon \lambda_{hj} +\varepsilon \gamma_{kl}$ in the case both or one tend to infinity. Therefore, $(\Sigma^{(r)},\Lambda^{(r)},\Gamma^{(r)})$ is not an infimizing sequence. Thus, (\ref{pb_GML}) is equivalent to minimize $\tilde \ell$ over the set 
 \al{\Cc_1=\{ (\Sigma,\Lambda,\Gamma) &\hbox{ s.t. } \Sigma\in\Qc_{m,n}^+, \nn\\ & \lambda_{hj}\in[0,\lambda_M], \; \gamma_{kl}\in[0,\gamma_M]\}\nn}
where $\lambda_M,\gamma_M>0$ are some constants sufficiently large.\\
$\bullet$ {\em Case 2:} $(\Sigma^{(r)},\Lambda^{(r)},\Gamma^{(r)})\in \Cc_1$ is a convergent sequence such that there exists at least one   $(\bar h ,\bar k,\bar j,\bar l)\in\Tc$ such that $\lambda^{(r)}_{\bar h\bar j }\rightarrow 0$ and $\gamma^{(r)}_{\bar k\bar l }\rightarrow 0$. Then, we have
\al{\lim_{r\rightarrow \infty}& \tilde \ell(y^N ,\Sigma^{(r)},\Lambda^{(r)},\Gamma^{(r)})\nn \\ &\geq \lim_{r\rightarrow \infty}  \ell(y^N;\Sigma^{(r)})-\hspace{-0.4cm}\sum_{(h,k,j,l)\in\Tc} \hspace{-0.4cm}\alpha_{hk,jl}\log \max\{\lambda_{hj}^{(r)},\gamma_{kl}^{(r)}\} \nn\\
&\geq c -  \lim_{r\rightarrow \infty} \sum_{(h,k,j,l)\in\Tc} \hspace{-0.4cm} \alpha_{hk,jl}\log\max\{\lambda_{hj}^{(r)},\gamma_{kl}^{(r)}\}=\infty \nn}
where we have exploited the following facts: $g(\Sigma;\Lambda,\Gamma)\geq 0$, $\varepsilon \lambda_{hj}$ and $\varepsilon \gamma_{kl}$ are nonnegative; $\ell(y^N;\Sigma)\geq c $ for any $\Sigma\in\bar \Qc_{m,n}^+$; $\alpha_{hk,jl}>0$; $\log\max\{\lambda_{\bar h \bar j}^{(r)},\gamma_{\bar k \bar l}^{(r)}\}\rightarrow -\infty$  . Accordingly, $(\Sigma^{(r)},\Lambda^{(r)},\Gamma^{(r)})$ is not an infimizing sequence. Thus, (\ref{pb_GML}) is equivalent to  minimize $\tilde \ell$ over one of the following sets: 
\al{\label{set_C2_1}\Cc_2&=\{ (\Sigma,  \Lambda,  \Gamma)  \hbox{ s.t. } \Sigma\in \Qc_{m,n}^+, \nn\\ &\hspace{2.5cm} \lambda_{hj}\in[\beta,\lambda_M],\;\gamma_{kl}\in[0,\gamma_M]\} \\ \label{set_C2_2} \Cc_2&=\{ (\Sigma,  \Lambda,  \Gamma)  \hbox{ s.t. } \Sigma\in \Qc_{m,n}^+, \nn\\ &\hspace{2.5cm}  \lambda_{hj}\in[0,\lambda_M],\;\gamma_{kl}\in[\beta,\gamma_M]\}}
where $\beta>0$ is a sufficiently small constant. Notice that $\max\{\lambda_{hj},\gamma_{kl}\}> 0$ for any $\Lambda$ and $\Gamma$ in $\Cc_2$. Without loss of generality we consider the set in (\ref{set_C2_1}). \\
$\bullet$ {\em Case 3:} $(\Sigma^{(r)},\Lambda^{(r)},\Gamma^{(r)})\in \Cc_2$ is a convergent sequence such that $\|\Sigma^{(r)}\|\rightarrow \infty$. Then, it is not difficult to see that $\ell(y^N; \Sigma^{(r)})\rightarrow \infty$, see Lemma 5.4  in \cite{FERRANTE_TIME_AND_SPECTRAL_2012}. Accordingly, we have 
\al{&\lim_{r\rightarrow \infty}  \tilde \ell(y^N ,\Sigma^{(r)},\Lambda^{(r)},\Gamma^{(r)})\nn \\ &\geq \lim_{r\rightarrow \infty}  \ell(y^N;\Sigma^{(r)})-\hspace{-0.4cm}\sum_{(h,k,j,l)\in\Tc}\hspace{-0.4cm} \alpha_{hk,jl}\log \max\{\lambda_{hj}^{(r)},\gamma_{kl}^{(r)}\} =\infty \nn}
where we have exploited the facts that $g(\Sigma;\Lambda,\Gamma)\geq 0$, $\varepsilon \lambda_{hj}\geq 0$, $\varepsilon \gamma_{kl}\geq 0$ and the terms $\alpha_{hk,jl}\log\max\{\lambda_{hj}^{(r)},\gamma_{kl}^{(r)}\}$ take finite values in $\Cc_2$. Accordingly, $(\Sigma^{(r)},\Lambda^{(r)},\Gamma^{(r)})$ is not an infimizing sequence. Thus, (\ref{pb_GML}) is equivalent to  minimize $\tilde \ell$ over the set: 
\al{\Cc_3=\{ (\Sigma,  \Lambda,  \Gamma)& \hbox{ s.t. } \Sigma\in \Qc_{m,n}^+, \Sigma \preceq \mu I\nn\\ & \lambda_{hj}\in[\beta,\lambda_M],\;\gamma_{kl}\in[0,\gamma_M]\}} where $\mu>0$ is a sufficiently large constant.\\
$\bullet$ {\em Case 4:} $(\Sigma^{(r)},\Lambda^{(r)},\Gamma^{(r)})\in \Cc_3$ converges to $(\Sigma_B,\Lambda_B,\Gamma_B)$ where $\Sigma_B$ lies on the boundary of $\Qc_{m,n}^+$, $\lambda_{B,hj}\in[\beta,\lambda_M]$ and $\gamma_{B,kl}\in [0,\gamma_M]$. First, note that the terms $g(\Sigma^{(r)};\Lambda^{(r)},\Gamma^{(r)})$, $\alpha_{hk,jl}\log \max\{\lambda_{hj}^{(r)},\gamma_{kl}^{(r)}\}$ and $\varepsilon (\lambda_{hj}^{(r)}+ \gamma_{kl}^{(r)})$ tend to finite values as  
$(\Sigma^{(r)},\Lambda^{(r)},\Gamma^{(r)})\rightarrow (\Sigma_B,\Lambda_B,\Gamma_B)$. Therefore, we just need to analyze the behaviour of  $\ell(y^N; \Sigma^{(r)})$. If $|\Sigma_B(e^{i\vartheta})|=0$ for any $\vartheta\in[-\pi,\pi]$, then $\ell(y^N; \Sigma^{(r)})\rightarrow \infty$, see Lemma 5.3 in \cite{FERRANTE_TIME_AND_SPECTRAL_2012}; thus  $\tilde \ell(y^N,\Sigma ^{(r)},\Lambda^{(r)},\Gamma^{(r)})\rightarrow \infty$, i.e. it is not an infimizing sequence. In the case that $|\Sigma_B(e^{i\vartheta})|$ is different from the null function we have $\ell(y^N; \Sigma ^{(r)})$, and thus also $\tilde \ell(y^N; \Sigma ^{(r)},\Lambda^{(r)},\Gamma^{(r)})$, converges to a finite value, see again Lemma 5.3 in \cite{FERRANTE_TIME_AND_SPECTRAL_2012};. On the other hand, the first variation of $\tilde \ell$ with respect to $\Sigma$ in $\Sigma_B$ along the direction $I$ (i.e. towards the interior of $\bar \Cc$) is
\al{\lim_{\varepsilon \downarrow 0}&\frac{\tilde \ell(y^N, \Sigma_B+\varepsilon I;\Lambda_B ,\Gamma_B)-\tilde \ell(y^N, \Sigma_B;\Lambda_B,\Gamma_B)}{\varepsilon}\nn\\ &=\lim_{\varepsilon \downarrow 0}\frac{  \ell(y^N; \Sigma_B+\varepsilon I)- \ell(y^N; \Sigma_B)}{\varepsilon}\nn\\ & +\lim_{\varepsilon \downarrow 0}\frac{  g( \Sigma_B+\varepsilon I;\Lambda_B,\Gamma_B)- g (\Sigma_B;\Lambda_B,\Gamma_B)}{\varepsilon}=-\infty \nn} 
where we have exploited the fact that the first term in the summation tends to $-\infty$, see Theorem 5.2 in \cite{FERRANTE_TIME_AND_SPECTRAL_2012}, and the second one converges to a bounded value. Accordingly, $(\Sigma^{(r)},\Lambda^{(r)},\Gamma^{(r)})$ is not an infimizing sequence. Therefore, (\ref{pb_GML}) is equivalent to minimize $\tilde \ell$ over the set \al{\Cc^\star=\{ (\Sigma,  \Lambda,  \Gamma)& \hbox{ s.t. } \Sigma\in \Qc_{m,n}^+, \, \upsilon I \preceq \Sigma \preceq \mu I\nn\\ & \lambda_{hj}\in[\beta,\lambda_M],\;\gamma_{kl}\in[0,\gamma_M]\}} where $\upsilon>0$ is a sufficiently small constant.\\
Since $\tilde \ell$ is a continuous function over $\Cc^\star$ and the latter is a compact set, by the Weierstrass theorem we conclude that there exists a point of minimum.
\qed\\

It is worth noting that $\tilde \ell$ is not a convex function, therefore the computation of a point of minimum is difficult. On the other hand, we will see that it is possible to find a point of minimum of $\tilde \ell$ with respect to $\Sigma$, $\Lambda$ and $\Gamma$, separately. For this reason, we propose the following three-step sequential algorithm for finding a coordinatewise minimum of $\tilde \ell$ (see Definition \ref{def_C_MIN} below): 
\al{\label{3step_sigma}\hat \Sigma^{(r+1)} =& \underset{  \Sigma}{\;\mathrm{argmin}}\;  \tilde \ell ( y^N,\Sigma, \hat \Lambda^{(r)},\hat \Gamma^{(r)})\\ 
\label{3step_lambda}\hat  \Lambda^{(r+1)} =& \underset{\Lambda }{\;\mathrm{argmin}}\; \tilde \ell ( y^N,\hat \Sigma^{(r+1)},\Lambda,\hat \Gamma^{(r)})\\
\label{3step_gamma}\hat{\Gamma}^{(r+1)} =& \underset{\Gamma }{\;\mathrm{argmin}}\; \tilde \ell ( y^N,\hat \Sigma^{(r+1)},\hat \Lambda^{(r+1)},\Gamma).}
Step (\ref{3step_sigma}) is the MAP estimator of $\Sigma$ given the current choice of $\Lambda,\Gamma$ which is equivalent to (\ref{reg_ml}). By Proposition \ref{prop_reg_MLA}, the latter always admits a unique solution which can be computed by means of a projective gradient algorithm \citep{REWEIGHTED}. 
Step (\ref{3step_lambda}) is the estimator of $\Lambda$ using the current MAP estimate of $\Sigma$ and the current choice of $\Gamma$, while step (\ref{3step_gamma}) is the estimator of $\Gamma$ using the current MAP estimate of $\Sigma$ and the current choice of $\Lambda$.
It is not difficult to see that the optimization of $\Lambda$ can be made independently for each entry, so we have 
{ \al{ \label{max_1} \hat{\lambda}^{(r+1)}_{hj}  = & \underset{ \lambda_{hj}\geq 0}{\mathrm{argmin}} \sum_{k\geq  l}  \max\{\lambda_{hj},\gamma_{kl}\}q_{hk,jl}( \Sigma)\nn\\ &-\alpha_{hk,jl} \log \max\{\lambda_{hj},\gamma_{kl}\} +\varepsilon \lambda_{hj} }
} where $\Sigma=\hat\Sigma^{(r+1)}$ and $\gamma_{kl}=\hat \gamma_{kl}^{(r)}$.
In a similar way, the optimization in (\ref{3step_gamma}) is equivalent to perform optimization for each entry independently:
{  \al{ \label{max_2} \hat{\gamma}^{(r+1)}_{kl} = &  \underset{\gamma_{kl}\geq 0}{\mathrm{argmin}} \sum_{h\geq j}\max\{  \gamma_{kl},\lambda_{hj}\} q_{hk,jl}( \Sigma )\nn\\ &-\alpha_{hk,jl} \log \max\{\lambda_{hj},\gamma_{kl}\}+\varepsilon \gamma_{kl} }
}  where $\Sigma=\hat\Sigma^{(r+1)}$ and $\lambda_{hj}=\hat \lambda_{hj}^{(r+1)}$.
 
\pro \label{prop_est_lg} Let $\tilde \gamma_u$, with $u=1\ldots \tilde u_2 $ and $\tilde u_2\leq m_2(m_2+1)/2$, be a reordering of the weights $\gamma_{kl}$, with $k\geq l$, such that $\tilde \gamma_{u+1}>\tilde \gamma_u$, i.e. if there are more than one weight taking the same value then the latter is taken only once in the reordering. Define the sets
\al{\Cc_u&=\{ (k,l) \hbox{ s.t. } k\geq l, \; \gamma_{kl}\leq \tilde \gamma_u \}\nn\\
\Nc_1&=\{ \tilde \gamma_u \hbox{ s.t. } u=1\ldots\tilde u_2 \}\nn\\
\Nc_2&=\{ \tilde \lambda_{hj,u} \hbox{ s.t. } u=0,1\ldots\tilde u_2 \}\nn } where 
\al{\tilde \lambda_{hj,u}= \left\{\begin{array}{ll}0 ,& u=0  \\ \frac{\sum_{(k,l)\in\Cc_u} \alpha_{hk,jl}}{\sum_{(k,l)\in\Cc_u} q_{hk,jl} (\Sigma)+\varepsilon}, & u>0. \end{array}\right.} Then, Problem (\ref{max_1}) admits solution and all the points of minimum are in the set 
\al{\Mc_1=\Nc_1 \cup\Nc_2.\nn}
\epro
 
\proof The objective function in (\ref{max_1}) is continuous but not necessarily differentiable in the points in $\Nc_1$. More precisely, the latter is convex and differentiable over the intervals $(0,\tilde \gamma_{1} )$, $(\tilde \gamma_u,\tilde \gamma_{u+1} )$ with $u=1\ldots \tilde u_2-1$, and $(\tilde \gamma_{\tilde u_2},\infty) $. For $\lambda_{hj}\in [\tilde \gamma_u,\tilde \gamma_{u+1} ]$, with $u=1\ldots \tilde u_2-1$, we have 
\al{\min\{\lambda_{hj},\gamma_{kl}\}= \left\{\begin{array}{ll}\lambda_{hj}, & (k,l)\in\Cc_u \\ \gamma_{kl}, & (k,l)\notin \Cc_u \end{array}\right.} and the objective function takes the form
\al{ \sum_{(k,l)\in\Cc_u} \lambda_{hj} q_{hk,jl}(\Sigma) -\alpha_{hk,jl} \log\lambda_{hj} +\varepsilon \lambda_{hj}+c } where $c$ is a constant not depending on $\lambda_{hj}$. The latter is strictly convex and its unique stationary point, if it {\mz exists}, is given by setting its first derivative equal to zero:
\al{ \sum_{(k,l)\in\Cc_u}  q_{hk,jl}(\Sigma) -\alpha_{hk,jl}\lambda_{hj}^{-1} +\varepsilon =0} and the solution of the above equality is $\tilde \lambda_{hj,u}$. If $\tilde \lambda_{hj,u}\in[\tilde \gamma_u,\tilde \gamma_{u+1}]$, then the  stationary point {\mz exists} and it is also the unique minimum of the objective function (\ref{max_1}) over   $[\tilde \gamma_u,\tilde \gamma_{u+1}]$. Otherwise, the point of minimum is $\tilde \gamma_u$ or $\tilde \gamma_{u+1}$. Therefore, the point of minimum of the objective function in (\ref{max_1}) over   $[\tilde \gamma_u,\tilde \gamma_{u+1}]$ is in the set $\{\tilde \lambda_{hj,u},\tilde \gamma_u, \tilde \gamma_{u+1}\}$. If $\lambda_{hj}\in[\tilde\gamma_{\tilde u_2},\infty)$ using similar arguments of before and the fact that the function tends to infinity as $\lambda_{hj}\rightarrow \infty$, we conclude that the point of minimum of the objective function in (\ref{max_1}) over  $[\tilde \gamma_{\tilde u_2},\infty)$ is in the set $\{\tilde \lambda_{hj,\tilde u_2},\tilde \gamma_{\tilde u_2}\}$. In the case that $\tilde \gamma_{1}>0$ we also have the interval $[0,\tilde \gamma_1]$. The objective function over this interval is equal to $\varepsilon \lambda_{hj} +c$ where $c$ is a constant not depending on $\lambda_{hj}$. Therefore, the point of minimum of the objective function in (\ref{max_1}) over $[0,\tilde \gamma_1]$ is $\tilde \lambda_{hj,0}=0$. We conclude that all the points of minimum of the objective function in (\ref{max_1}) are in the set $\Mc_1$.  
\qed\\ 

Proposition \ref{prop_est_lg} also provides  a simple way to solve (\ref{max_1}): it is just required to evaluate the objective function for the $\tilde u_2$ points in $\Mc_1$ and choose  one minimizing the function. 
 
\pro \label{prop_est_lg2} Let $\tilde \lambda_u$, with $u=1\ldots \tilde u_1 $ and $\tilde u_1\leq m_1(m_1+1)/2$, be a reordering of the weights $\lambda_{hj}$, with $h\geq j$, such that $\tilde \lambda_{u+1}>\tilde \lambda_u$, i.e. if there are more than one weight taking the same value then the latter is taken only once in the reordering. Define the sets
\al{\Cc_u&=\{ (h,j) \hbox{ s.t. } h\geq j, \; \lambda_{hj}\leq \tilde \lambda_u \}\nn\\
\Nc_1&=\{ \tilde \lambda_u \hbox{ s.t. } u=1\ldots\tilde u_1 \}\nn\\
\Nc_2&=\{ \tilde \gamma_{kl,u} \hbox{ s.t. } u=0,1\ldots\tilde u_1 \}\nn } where 
\al{\tilde \gamma_{kl,u}= \left\{\begin{array}{ll}0 ,& u=0  \\ \frac{\sum_{(h,j)\in\Cc_u} \alpha_{hk,jl}}{\sum_{(h,j)\in\Cc_u} q_{hk,jl} (\Sigma)+\varepsilon}, & u>0. \end{array}\right.} Then, Problem (\ref{max_2}) admits solution and all the points of minimum are in the set 
\al{\Mc_2=\Nc_1 \cup\Nc_2.\nn}
\epro

\proof The proof is similar to one of Proposition \ref{prop_est_lg}. \qed\\

As before, Proposition \ref{prop_est_lg2} provides a simple way to solve (\ref{max_2}). The three-step sequential procedure is summarized in Algorithm \ref{algo:RWS} where the iterative scheme ends when $\tilde \ell $ does not significantly change according to a fixed tolerance $\epsilon>0$.
 \begin{algorithm}
\caption{ }
\label{algo:RWS}
\begin{algorithmic}[1] \small
\STATE $r=0$
\STATE Initialize $\hat \lambda_{hj}^{(0)}$ and $\hat \gamma_{kl}^{(0)}$ 
\REPEAT 
\STATE Compute $\hat \Sigma^{(r+1)}$ by solving Problem (\ref{reg_ml}) with $\lambda_{hj}=\hat{\lambda}_{hj}^{(r)}$ and $\gamma_{kl}=\hat{\gamma}_{kl}^{(r)}$
\STATE Compute $q_{hk,jl}(\hat \Sigma^{(r+1)})$ 
\STATE Compute the weights $\hat \lambda_{hj}^{(r+1)}$ by solving Problem (\ref{max_1}) with $\Gamma=\hat \Gamma^{(r)}$ 
\STATE  Compute the weights $\hat \gamma_{kl}^{(r+1)}$  by solving Problem (\ref{max_2}) with $\Lambda=\hat\Lambda^{(r+1)}$ 
\STATE Compute $\tilde \ell^{(r+1)}:=\tilde \ell (y^N,\hat \Sigma^{(r+1)},\hat \Lambda^{(r+1)},\hat\Gamma^{(r+1)})$
\STATE $r\leftarrow r+1$
\UNTIL{$|\tilde \ell^{(r)}-\tilde \ell^{(r-1)}|\leq \epsilon $  }
\STATE $\hat \Phi=(\hat \Sigma^{(r+1)})^{-1}$ 
\end{algorithmic}
\end{algorithm}

A variation of the previous algorithm is to modify the sequence of (\ref{3step_sigma})-(\ref{3step_gamma}):
\al{\label{3step_sigmab}\hat \Sigma^{(r+1)} =& \underset{  \Sigma}{\;\mathrm{argmin}}\;  \tilde \ell ( y^N,\Sigma, \hat \Lambda^{(r)},\hat \Gamma^{(r)})\\ 
\label{3step_gammab}\hat{\Gamma}^{(r+1)} =& \underset{\Gamma }{\;\mathrm{argmin}}\; \tilde \ell ( y^N,\hat \Sigma^{(r+1)},\hat \Lambda^{(r)},\Gamma)\\
\label{3step_lambdab}\hat  \Lambda^{(r+1)} =& \underset{\Lambda }{\;\mathrm{argmin}}\; \tilde \ell ( y^N,\hat \Sigma^{(r+1)},\Lambda,\hat \Gamma^{(r+1)})}
leading to Algorithm \ref{algo:RWS2}. {\mz Notice that there is no essential difference between the two algorithms: only the order of optimizing $\Lambda$ and $\Gamma$ has been changed.}
 
\begin{algorithm}
\caption{ }
\label{algo:RWS2}
\begin{algorithmic}[1] \small
\STATE $r=0$
\STATE Initialize $\hat \lambda_{hj}^{(0)}$ and $\hat \gamma_{kl}^{(0)}$ 
\REPEAT 
\STATE Compute $\hat \Sigma^{(r+1)}$ by solving Problem (\ref{reg_ml}) with $\lambda_{hj}=\hat{\lambda}_{hj}^{(r)}$ and $\gamma_{kl}=\hat{\gamma}_{kl}^{(r)}$
\STATE Compute $q_{hk,jl}(\hat \Sigma^{(r+1)})$
\STATE  Compute the weights $\hat \gamma_{kl}^{(r+1)}$   by solving Problem (\ref{max_2}) with $\Lambda=\hat\Lambda^{(r)}$ 
\STATE Compute the weights $\hat \lambda_{hj}^{(r+1)}$  by solving Problem (\ref{max_1}) with $\Gamma=\hat \Gamma^{(r+1)}$ 
\STATE Compute $\tilde \ell^{(r+1)}:=\tilde \ell (y^N,\hat \Sigma^{(r+1)},\hat \Lambda^{(r+1)},\hat \Gamma^{(r+1)})$ 
\STATE $r\leftarrow r+1$
\UNTIL{$|\tilde \ell^{(r)}-\tilde \ell^{(r-1)}|\leq \epsilon $  }
\STATE $\hat \Phi=(\hat \Sigma^{(r+1)})^{-1}$ 
\end{algorithmic}
\end{algorithm}

Next we show the properties of the limit points of the proposed algorithms. 

\df {\mz\citep{SEI}} \label{def_C_MIN}We say that $(\bar \Sigma, \bar \Lambda,\bar \Gamma)\in\Cc$ is a coordinatewise minimum point of $\tilde \ell$ if the following conditions hold:
\al{\label{cond_w_min}&\tilde \ell (y^N , \bar \Sigma; \bar \Lambda,\bar \Gamma)\leq \tilde \ell (y^N , \bar \Sigma+\Xi; \bar \Lambda,\bar \Gamma)\nn\\
&\tilde \ell (y^N , \bar \Sigma; \bar \Lambda,\bar \Gamma)\leq \tilde \ell (y^N , \bar \Sigma; \bar \Lambda+\Delta,\bar \Gamma)\nn\\
&\tilde \ell (y^N , \bar \Sigma; \bar \Lambda,\bar \Gamma)\leq \tilde \ell (y^N , \bar \Sigma; \bar \Lambda,\bar \Gamma+\Upsilon)} for any $\Xi, \Delta,\Upsilon$ such that $\bar\Sigma +\Xi\in\Qc_{m,n}^+$, 
$\lambda_{hj}+\delta_{hj}\in[0,\lambda_M]$, $\gamma_{kl}+\upsilon_{kl}\in[0,\gamma_M]$ and $\delta_{hj}$, $\upsilon_{kl}$ are the entries in position $(h,j)$ and $(k,l)$ of $\Delta$ and $\Upsilon$, respectively. 
\edf

Clearly a coordinatewise minimum point is also a stationary point, but the converse is not true, i.e. coordinatewise minimum point is a stronger property than  stationary point.

\corr \label{prop_limit_point}Let $(\hat \Sigma^{(r)},\hat \Lambda^{(r)},\hat\Gamma ^{(r)})$, with $n\in\Ns$, be the sequence generated by one of the previous algorithms. {\mz If there is a limit point of such a sequence, then it  is a coordinatewise minimum point of $\tilde \ell$.} \ecorr
\proof The statement follows from the fact that at each step of the sequential procedure we find a point of minimum for $\tilde \ell$ with respect to one variable.\qed\\

In the sequential steps  (\ref{3step_sigma})-(\ref{3step_gamma}) and 
 (\ref{3step_sigmab})-(\ref{3step_lambdab}) we have to select the initial conditions for $\Lambda$ and $\Gamma$ that is $\Lambda^{(0)}$ and $\Gamma^{(0)}$. Clearly, the better the initialization is, the better the final estimate of $\Sigma$ will be. The idea is to estimate  $\Lambda^{(0)}$ and $\Gamma^{(0)}$ from the a preliminary estimate of $\Sigma$. More precisely, let $\hat \Phi_{ME}=(\hat \Sigma_{B})^{-1}$ denote {\mz the estimator of $\Phi$ obtained by solving (\ref{reg_ml}) without regularization, i.e., with $\Lambda=0$ and $\Gamma=0$, cf. Section \ref{sec:ME}}. Then, we solve iteratively   
 \al{\hat  \Lambda_{\mathrm{init}}^{(q+1)} =& \underset{\Lambda }{\;\mathrm{argmin}}\; \tilde \ell ( y^N,\hat \Sigma_{B},\Lambda,\hat \Gamma_{\mathrm{init}}^{(q)})\nn\\
\hat{\Gamma}_{\mathrm{init}}^{(q+1)} =& \underset{\Gamma }{\;\mathrm{argmin}}\; \tilde \ell ( y^N,\hat \Sigma_{B},\hat \Lambda_{\mathrm{init}}^{(q+1)},\Gamma)\nn}
where $\hat \Gamma_{\mathrm{init}}^{(0)}$ is a matrix of ones and the iterative procedure stops when $\| \hat \Lambda_{\mathrm{init}}^{(q+1)}-\hat \Lambda_{\mathrm{init}}^{(q)} \|\leq \tilde \epsilon$ and $\| \hat \Gamma_{\mathrm{init}}^{(q+1)}-\hat \Gamma_{\mathrm{init}}^{(q)} \|\leq \tilde \epsilon$ for some $\tilde \epsilon>0$ sufficiently small. Finally, we set $\Lambda^{(0)}=\hat \Lambda_{\mathrm{init}}^{(q)}$ and $\hat \Gamma^{(0)}= \hat \Gamma_{\mathrm{init}}^{(q)}$.

\section{Multiplicative prior}\label{sec:additive2}
We consider the possibility to use a penalty function which generalizes the one in \cite{CDC_KRON} for learning static KGM. More precisely, the idea is to replace in (\ref{hA}) the term $\max\{\lambda_{hj},\gamma_{kl}\}$ with the multiplicative term $\lambda_{hj}\gamma_{kl}$. Therefore, we consider the penalty term:
  \al{\label{hB}   g &(\Sigma; \Lambda,\Gamma) =\hspace{-0.2cm}\sum_{(h,k,j,l)\in\Tc}  \hspace{-0.4cm} \lambda_{hj}\gamma_{kl}q_{hk,jl}(\Sigma).}  
Such a penalty term is in the same spirit of the prior proposed in \cite{bonilla2008multi} and \cite{yu2009large} for multi-task learning and collaborative filtering, respectively. In what follows, we consider the regularized ML problem in (\ref{reg_ml}) with $g(\Sigma;\Lambda,\Gamma)$ defined in (\ref{hB}). Taking the Bayesian perspective of Section \ref{sec:additive}, we model $\Sigma$ as a stochastic process taking values in $\Qc_{m,n}^+$ and such that (\ref{form_indip}) holds. We define the ``multiplicative prior''
\al{ p((\Sigma)_{hk,jl}|\lambda_{hj},\gamma_{kl})=\frac{e^{-\lambda_{hj}\gamma_{kl}q_{hk,jl}(\Sigma)}}{c_{hk,jl}} \nn}
where $c_{hk,jl}>0$ is the normalizing constant.   Also in this case we model $\Lambda$ and $\Gamma$ as random matrices whose PDF is given by (\ref{form_LG_1})-(\ref{form_LG_3}). Then, an upper bound for $\ell(y^N,\Sigma,\Lambda,\Gamma)$ is 
\al{\label{upper_bound_sommaB} \tilde \ell &(    y^N,  \Sigma, \Lambda,\Gamma) := \ell(y^N; \Sigma) +g(\Sigma;\Lambda,\Gamma)\nn\\&-\hspace{-0.4cm} \sum_{(h,k,j,l)\in\Tc} \hspace{-0.4cm} \alpha_{hk,jl} \log ( \lambda_{hj} \gamma_{kl} ) +\varepsilon\sum_{h\geq j} \lambda_{hj}+\varepsilon\sum_{k\geq l} \gamma_{kl}. }
It is not difficult to see that the problem in (\ref{pb_GML}) with (\ref{upper_bound_sommaB}) {\mz admits} solution. Then, the optimization problem can be solved by the sequential procedures (\ref{3step_sigma})-(\ref{3step_gamma}) and
 (\ref{3step_sigmab})-(\ref{3step_lambdab}). Also in this case the optimization of $\Lambda$ and $\Gamma$ can be made independently for each entry: 
 \al{ \label{max_1b} \hat{\lambda}^{(r+1)}_{hj} & =  \underset{ \lambda_{hj}\geq 0}{\mathrm{argmin}} \sum_{k\geq  l}   \lambda_{hj} \gamma_{kl} q_{hk,jl}( \Sigma)\nn\\ 
 &\hspace{1cm}-\alpha_{hk,jl} \log (\lambda_{hj}\gamma_{kl}) +\varepsilon \lambda_{hj}  \\
  \label{max_2b}   \hat{\gamma}^{(r+1)}_{kl} &=    \underset{\gamma_{kl}\geq 0}{\mathrm{argmin}}  \sum_{h\geq j}   \gamma_{kl} \lambda_{hj} q_{hk,jl}( \Sigma )\nn\\ &\hspace{1cm}-\alpha_{hk,jl} \log (\lambda_{hj}\gamma_{kl})+\varepsilon \gamma_{kl}  
   } where $\Sigma=\hat\Sigma^{(r+1)}$, $\gamma_{kl}=\hat \gamma_{kl}^{(r)}$ in (\ref{max_1b}) and $\lambda_{hj}=\hat \lambda_{hj}^{(r+1)}$ in (\ref{max_2b}) if we consider the sequential scheme (\ref{3step_sigma})-(\ref{3step_gamma}). It is not difficult to prove that (\ref{max_1b}) and (\ref{max_2b}) admit  unique solution whose analytic expression are, respectively: 
  \al{ \hat \lambda_{hj}^{(r+1)} &=\left\{\begin{array}{ll}\frac{1}{2} \frac{m_2+m_2^2(2n+1)}{ \sum_{k\geq l}    \gamma_{kl} q_{hk,hl}( \Sigma) +\varepsilon} , & \hbox{ if $h=j$} \nn\\  \frac{m_2^2(2n+1)}{\sum_{k\geq  l}   \gamma_{kl} q_{hk,hl}(  \Sigma )+\varepsilon} , & \hspace{-0.2cm}\hbox{otherwise}\end{array}\right. \nn\\
\hat \gamma_{kl}^{(r+1)} &=\left\{\begin{array}{ll} \frac{1}{2}\frac{m_1+m_1^2(2n+1)}{ \sum_{h\geq j}   \lambda_{hj}  q_{hk,jl}(  \Sigma )+\varepsilon} , & \hbox{if $k=l$} \nn\\ \frac{m_1^2(2n+1)}{\sum_{h\geq j}  \lambda_{hj}   q_{hk,jl}(  \Sigma )+\varepsilon} , & \hbox{otherwise.}\end{array}\right. }  

\rema It is worth noting that the resulting sequential procedure is similar to an iterative reweighting scheme \citep{WIPF_2010,scheinberg2010sparse}: in Step 4 we compute the regularized ML estimator where (\ref{cond_sup}) is induced by the weighted penalty term $g(\Sigma;\Lambda,\Gamma)$; in Step 6 the weight  $\lambda_{hj}$ is inversely proportional to $\sum_{k\geq l} \gamma_{kl}q_{hk,jl}(\Sigma )$ which is the weighted $\ell_1$ norm of the $m_2\times m_2$ matrix block $(\Sigma)_{hk,jl}$ with $k,l\in\Vc_2$ of the current estimate of $\Sigma$;  in Step 7 the weight  $\gamma_{kl}$ is inversely proportional to $\sum_{h\geq j} \lambda_{hj}q_{hk,jl}(\Sigma )$ which is the weighted $\ell_1$ norm of the $m_1\times m_1$ matrix block $(\Sigma)_{hk,jl}$ with $h,j\in\Vc_1$ of the current estimate of $\Sigma$. {\mz The main difference is that in our method we have two priors to update sequentially while in \citep{WIPF_2010} only one prior is present. } \erema

The penalty function (\ref{hB}) is more appealing than the one in (\ref{hA}) because it guarantees the uniqueness of the minimum with respect to $\Lambda$ and $\Gamma$. However, as we will see in Section \ref{sec:sim}, the penalty (\ref{hB}) does not provide a good performance. Such an evidence can be justified as follows: {\mz assume data is generated} from a model whose PSD is such that condition (\ref{cond_sup}) holds with $(\bar h,\bar j)\notin \Ec_1$ and $(\bar k,\bar l)\in\Ec_2$. Then, the procedure will tend to assign a large weight corresponding to $q_{\bar hk,\bar jl}(\Sigma)$ for any $k\geq l$. Notice that the weight corresponding to $q_{\bar hk,\bar jl}(\Sigma)$ is $\lambda_{\bar h\bar j} \gamma_{kl}$. If the current  $\lambda_{\bar h \bar j}$ is not so large, then the optimization of $\gamma_{\bar k \bar l}$ leads to a value which is large in order to prune $q_{\bar h\bar k,\bar j\bar l}(\Sigma)$ to zero. Such a value for $\gamma_{\bar k \bar l}$ is wrong because it prunes to zero the nonnull entries $(\Sigma)_{h\bar k, j \bar l}$ with $(h,j)\neq (\bar h,\bar j)$.    
 
\section{A Maximum Entropy interpretation}\label{sec:ME}
We show that Problem (\ref{reg_ml}) with penalty term (\ref{hA}) or (\ref{hB}) is connected with a maximum entropy (ME) problem. We consider the Gaussian process $y$ of Section \ref{sec:pn_form} taking values in $\Rs^{m_1m_2}$. Given the data $y^N$, the Burg spectral estimator \citep{burg1975maximum} solves the following covariance extension problem:
\al{\label{pb_BURG} \hat \Phi_{ME}= &\underset{\Phi}{ \mathrm{argmax}} \,\frac{1}{4\pi}\int_{-\pi}^\pi \log|\Phi (e^{i\vartheta})| \mathrm d \vartheta\nn\\
& \hbox{ s.t. }  \Phi\succ 0\nn\\
&  \frac{1}{2\pi}\int_{-\pi}^\pi \Phi(e^{i\vartheta}) e^{-i\vartheta s} \mathrm d\vartheta =\hat R_s, \, \, s=0\ldots n.}
Such estimator is also known as ME estimator because the objective function is the differential entropy rate of the process with PSD $\Phi$. 
It is worth noting there exist also alternative objective functions for which a large body of literature has been produced, e.g. \cite{A_NEW_APPROACH_BYRNES_2000,Hellinger_Ferrante_Pavon,8359301,BETA,1603771}. In Problem (\ref{pb_BURG}) we impose that $\hat \Phi$ matches the first $n$ moments (i.e. covariance lags) $\hat R_s$'s which are estimated from $y^N$. Assume that we have some a priori knowledge about the reliability of $\hat R_s$'s. More precisely, we assume to know that: (i) the 
information between module $h$ and module $j$ (i.e. the covariances between the nodes in module $h$ and in module $j$, respectively) is not reliable if $(h,j)\notin\Ec_1$; (ii) the covariance  between node $k$ and node $l$ in each module is not reliable if $(k,l)\notin\Ec_2$. If we discard the non-reliable data, we obtain      
 \al{\label{pb_DEMP} & \hat \Phi_{ME} =  \underset{\Phi}{ \mathrm{argmax}} \,\frac{1}{4\pi}\int_{-\pi}^\pi \log|\Phi (e^{i\vartheta})| \mathrm d \vartheta\nn\\
& \hbox{ s.t. }  \Phi\succ 0\nn\\
&  \left(\frac{1}{2\pi}\int_{-\pi}^\pi \Phi(e^{i\vartheta}) e^{-i\vartheta s} \mathrm d\vartheta \right)_{hk,jl}=\left(\hat R_s\right)_{hk,jl}, \nn\\ & \;\; \;   \;\; \;   \;\; \;   \;\; \;   s=0\ldots n, \; \forall\, (h,j)\in \Ec_1 \hbox{ and } \forall \, (k,l)\in\Ec_2.}
Using the duality theory as in \cite{ARMA_GRAPH_AVVENTI}, it is not difficult to prove that the dual of (\ref{pb_DEMP}) is 
\al{ \label{pd_dual_stric}\hat \Sigma= &\underset{\Sigma\in \Qc_{m,n}}{ \mathrm{argmin}} \,\ell( y^N; \Sigma)\nn\\
& \hbox{ s.t. }  \Sigma\in\Qc_{m,n}^+\nn\\
& (\Sigma)_{hk,jl} =0\; \; \forall\, (h,j)\notin \Ec_1 \hbox{ and } \forall \, (k,l)\notin\Ec_2 } 
and $\hat \Phi_{ME}=(\hat \Sigma)^{-1}$. Therefore, Problem (\ref{pd_dual_stric}) searches an AR model of order $n$ satisfying  (\ref{cond_sup}). Notice that, the latter is imposed as hard constraint and  
 it is required to know in advance the topology of the KGM, i.e. $\Ec_1$ and $\Ec_2$. In Problem (\ref{reg_ml}), instead, (\ref{cond_sup}) is imposed as soft constraint and it is not required to know in advance the topology of the graph. {\mz Another important aspect is that we have a perfect partial covariance matching in (\ref{pb_DEMP})-(\ref{pd_dual_stric}), while we have an approximate covariance matching in (\ref{reg_ml}). The latter strategy has been successfully used in other spectral estimation problems, see  \cite{TRE,QUATTRO}.}
 
 \section{Simulation Results}\label{sec:sim}
 \subsection{Synthetic Data}
We compare the performance of the KGM estimators proposed in Section \ref{sec:additive}. {\mz The corresponding Matlab functions are available at \url{https://github.com/MattiaZ85/KR-AR-GM}.} We will use the following shorthand notations: S denotes the SGM estimator proposed in \cite{REWEIGHTED} with penalty term (\ref{hSPARSE}); K1 denotes Algorithm \ref{algo:RWS}; K2 denotes Algorithm \ref{algo:RWS2}. In all the aforementioned estimators we set $\varepsilon=10^{-3}$ and $\epsilon=10^{ -3}$. In what follows we consider four Monte Carlo studies constituted by $200$ experiments. In each experiment we generate randomly an AR stochastic process $y$ of dimension $m_1m_2=36$ and order $n=2$. {\mz The latter values are kept fixed in oder to obtain an homogeneous comparison among the Monte Carlo studies.} The PSD is denoted by $\Phi$ and is such that $\mathrm{supp}(\Phi^{-1})=E_1\otimes E_2$. We denote by $\eta_1$ and $\eta_2$ the fraction of ones in $E_1$ and $E_2$, respectively. These supports are chosen randomly at each run. Then, we generate a finite length sequence $y(1)\ldots y(N)$ with $N=1000$. To asses the performance of the estimators we compute the fraction of misspecified edges with respect to the true PSD:
{\mz \al{e_{SP}=\frac{\|E_1\otimes E_2-\hat{E}_1\otimes \hat{ E}_2\|_0}{m_1^2 m_2^2}}
where $\mathrm{supp}(\hat \Phi^{-1})=\hat E_1\otimes \hat E_2  $, $\hat \Phi$ denotes the estimator of $\Phi$ and $\|\cdot \|_0$ denotes the $\ell_0$ matrix norm.} Moreover, we compute the relative error of $\hat\Phi^{-1}$ with respect to the true inverse PSD $\Phi^{-1}$:
{\mz \al{err=\frac{\int_{-\pi}^{\pi} \|\hat \Phi^{-1}-\Phi^{-1}\|_F^2\mathrm d\vartheta}{\int_{-\pi}^{\pi} \|\Phi^{-1}\|^2_F\mathrm d\vartheta}} where $\|\cdot\|_F$ denotes the Frobenius norm. Notice that we have considered the aforementioned error rather than the one between $\hat\Phi$ with $\Phi$ because $\Phi^{-1}$ provides information more close to the graphical model: the off-diagonal entries of $\Phi^{-1}$ represents the ``weight functions'' corresponding to the edges of the graphical model.}
  \begin{figure}[htbp]
\centering
  \includegraphics[width=\columnwidth]{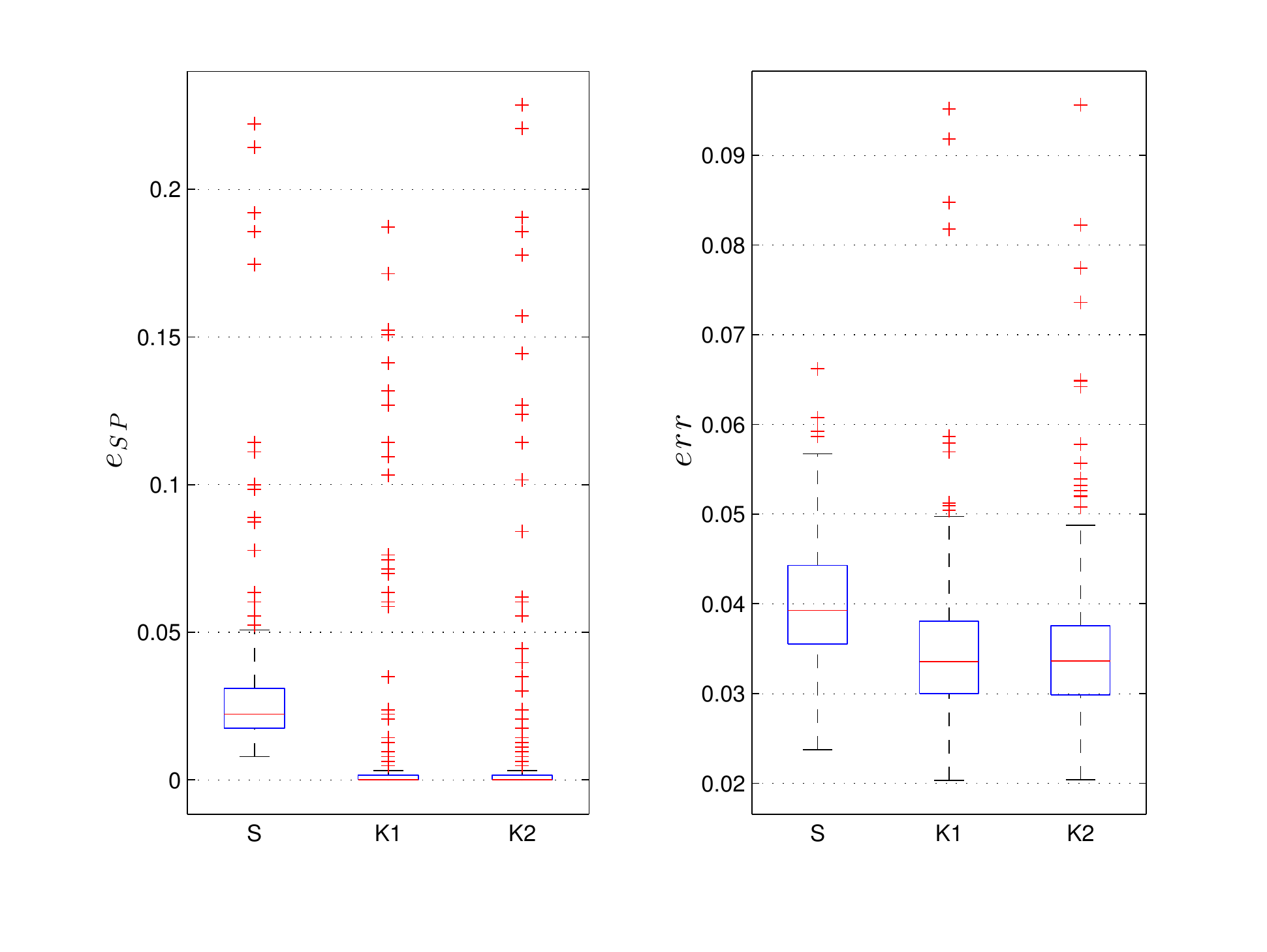}  
\caption{Boxplots of the fraction of misspecified edges (left) and the relative error (right) with $m_1=m_2=6$, $n=2$ and $\eta_1=\eta_2=0.3$.} \label{fig_due}
\end{figure}

\begin{figure}[htbp]
\centering
  \includegraphics[width=\columnwidth]{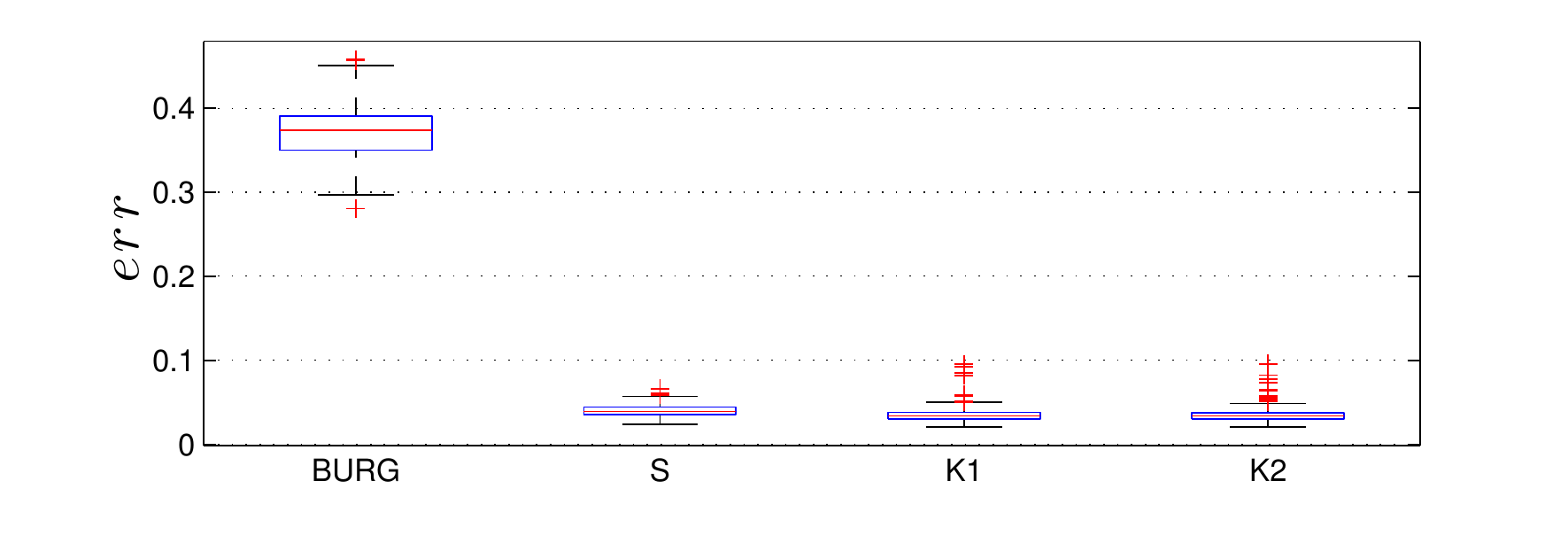}  
\caption{\mz Boxplots of the relative error of the Burg estimator (BURG), S, K1 and K2.} \label{fig_duebis}
\end{figure}

{\em First Monte Carlo study.} We set $m_1=m_2=6$ and $\eta_1=\eta_2=0.3$ that is $E_1$ and $E_2$ in the Kronecker decomposition have the same characteristics. Figure \ref{fig_due} shows the boxplots of the fraction of misspecified edges (left) and the relative error (right) for each estimator: the proposed estimators outperforms S. Moreover, K1 and K2 performs in the same way. This means that the sequence of the three optimization problems does not play any role. {\mz In Figure \ref{fig_duebis} we compare $err$ of the previous estimators with the one obtained with the Burg estimator, with $n=2$, which does not impose any kind of regularization for the network topology: the inferior performance of the latter when compared to S, K1 and K2 is more salient. This is because S, K1 and K2 search the optimal model over a suitable restricted model class and such a restriction depends on the regularizers. Finally, we have considered also the relative error between $\hat \Phi$ and $\Phi$: we have found a situation similar to the one corresponding to $err$. }

{\em Second Monte Carlo study.} We set $m_1=9$, $m_2=4$ and $\eta_1=\eta_2=0.3$. This means that the dimensions of $E_1$ and $E_2$ are different. \begin{figure}[htbp]
\centering
  \includegraphics[width=\columnwidth]{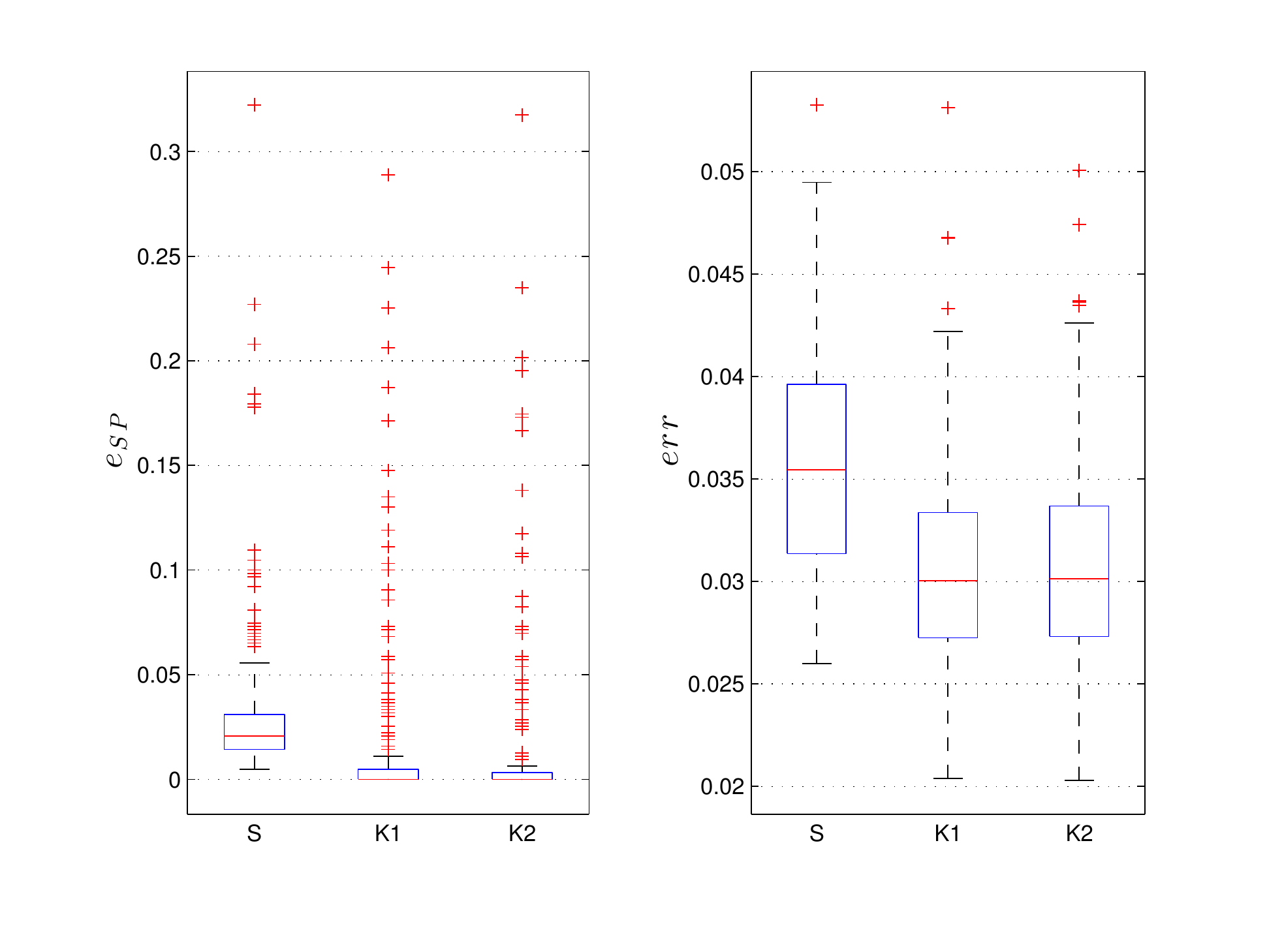}  
\caption{Boxplots of the fraction of misspecified edges (left) and the relative error (right) with $m_1=9$, $m_2=4$, $n=2$ and $\eta_1=\eta_2=0.3$.} \label{fig_otto}
\end{figure} Figure \ref{fig_otto} shows the boxplots of $e_{SP}$ (left) and $err$ (right)  for each estimator. Also in this case K1  and K2  outperform S. More precisely, the proposed estimators perform in the same way. Therefore, also in this case the sequence of the three optimization problems does not play a crucial role. We also have performed a Monte Carlo study with $m_1=4$ and $m_2=9$ (i.e. we have swapped $m_1$ and $m_2$) and we have obtained similar results.

\begin{figure}[htbp]
\centering
  \includegraphics[width=\columnwidth]{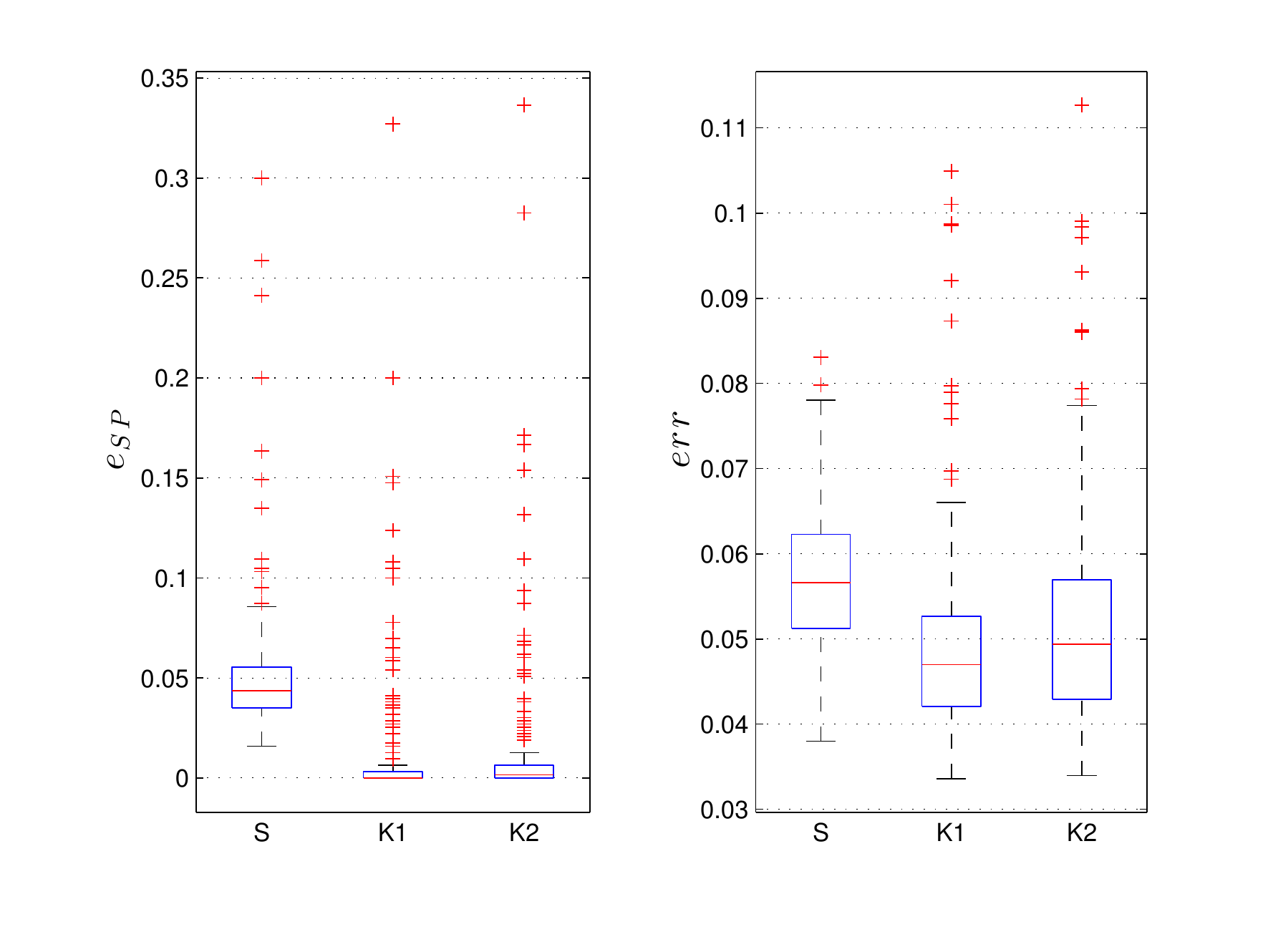}  
\caption{ Boxplots of the fraction of misspecified edges (left) and the relative error (right) with $m_1=m_2=6$, $n=2$, $\eta_1=0.3$ and $\eta_2=0.5$.} \label{fig_quattro}
\end{figure}

{\em Third Monte Carlo study.} We set $m_1=m_2=6$, $\eta_1=0.3$ and $\eta_2=0.5$ that is $E_1$ and $E_2$ differ from their degree of sparsity. In this situation $E_2$ is far from being sparse. 
Figure \ref{fig_quattro} 
shows the boxplots of $e_{SP}$ (left) and $err$ (right). K1 and K2  outperforms S. Moreover, K1 performs slightly better than  K2, however the performances are similar. Accordingly, the sequence in the optimization step does not play a crucial role. We also have performed a Monte Carlo study with $\eta_1=0.5$ and $\eta_2=0.3$ (i.e. we have swapped $\eta_1$ and $\eta_2$) and we have obtained a specular behavior: K2 performs slightly better than  K1 and the latter outperform S.

{\em Fourth Monte Carlo study}. \begin{figure}[htbp]
\centering
  \includegraphics[width=\columnwidth]{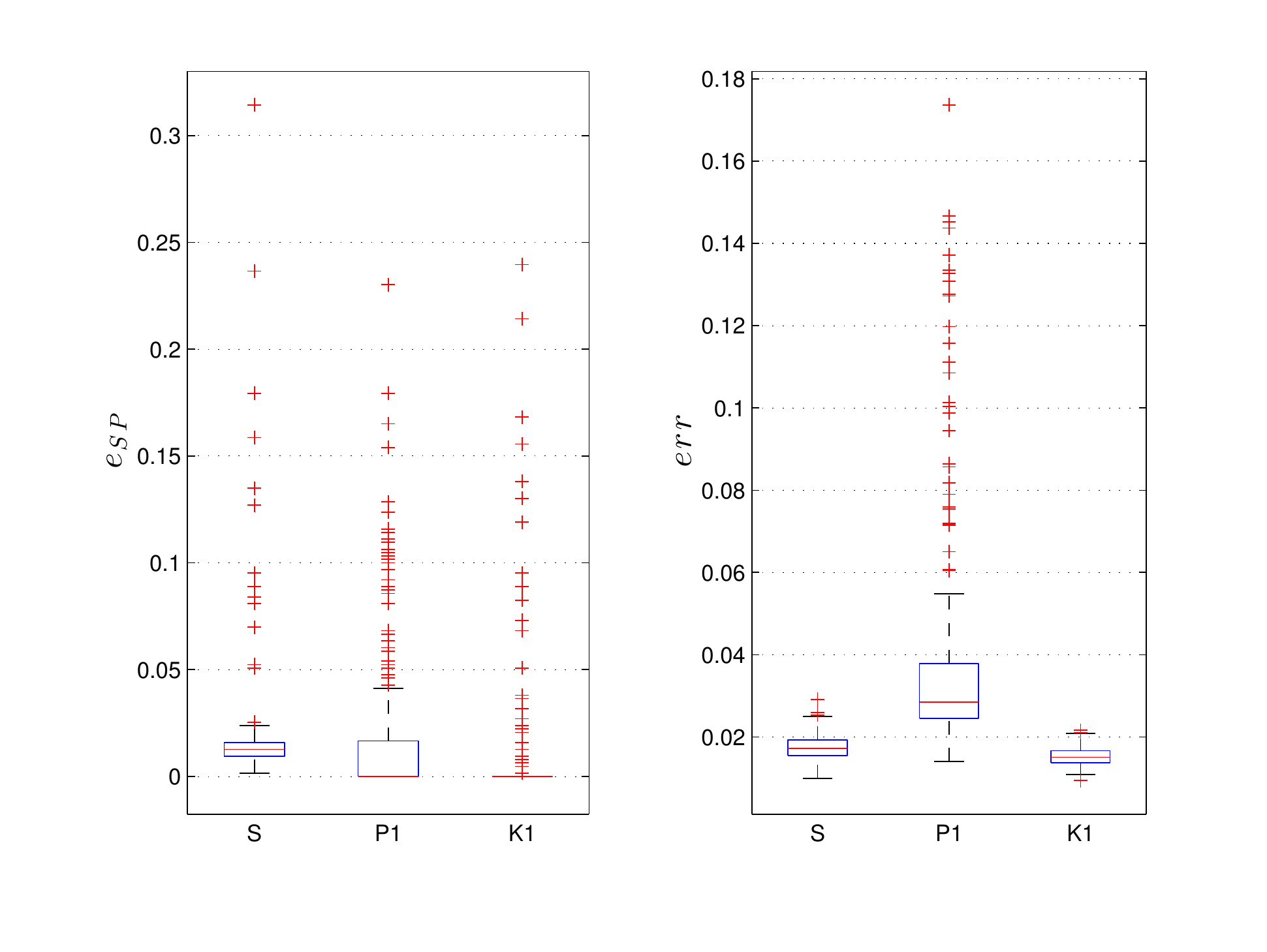}  
\caption{Boxplots of the fraction of misspecified edges (left) and the relative error (right) with $N=2000$, $m_1=m_2=6$, $n=2$, $\eta_1=0.3$ and $\eta_2=0.3$.} \label{fig_sei}
\end{figure} We compare the performance of S, K1 and P1. The latter denotes the  sequential procedure (\ref{3step_sigma})-(\ref{3step_gamma}), i.e. Algorithm \ref{algo:RWS}, which uses the multiplicative prior of Section \ref{sec:additive2}. Using the dataset of the first Monte Carlo study, we have found that P1 performs worse than S and K1. We have increased the size of the data $N=2000$, obtaining the results depicted in  Figure \ref{fig_sei}. K1 is still the best estimator. P1 outperforms S in terms of median of $e_{SP}$, but it is worse than S in terms of $err$. We obtained similar results using the sequential procedure of Algorithm \ref{algo:RWS2}. We conclude that the penalty term (\ref{hA}) is more effective than (\ref{hB}).

\begin{figure}[htbp]
\centering
  \includegraphics[width=\columnwidth]{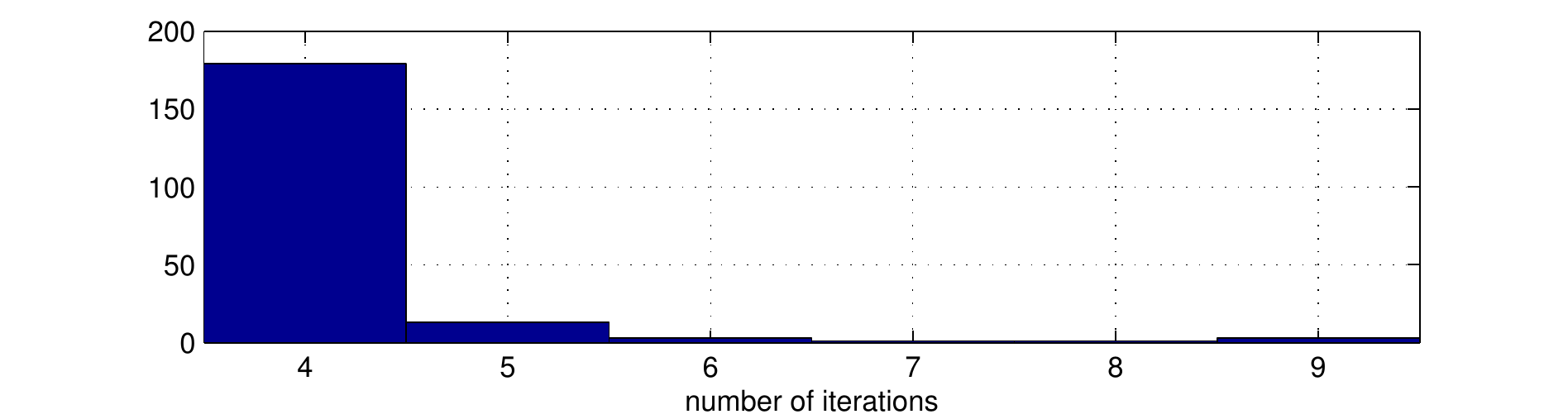}  
\caption{\mz Distribution of the number of performed iterations of the reweighting scheme for K1 in the first Monte Carlo study.} \label{fig_iter}
\end{figure}
{\mz {\em Computational complexity}. In the first Monte Carlo study the average computational time required to estimate the KGM using K1 is 243.67s, while the one using   K2 is 244.12s. These results have been obtained using a 6-Core AMD Opteron 8431 @2.4GHz processor. Similar results have been found in the other Monte Carlo studies. Moreover, Figure \ref{fig_iter} shows how the number of performed iterations of the reweighting scheme is distributed for K1. A similar distribution has been found for K2. In most of the cases, the number of required iterations is small, meaning that the convergence of the reweighting scheme is fast. The updating of $\Lambda$ and $\Gamma$ is very fast, indeed, as stated by Proposition \ref{prop_est_lg} and Proposition \ref{prop_est_lg2}, it has a closed form solution. The bottleneck is the computation of the solution to Problem (\ref{reg_ml}). Such solution can be found by considering a matrix optimization problem involving matrices of dimension $m(n+1)$, see \cite{REWEIGHTED} and \cite{SONGSIRI_TOP_SEL_2010}. The solution of the latter is found by a projected gradient whose computational complexity of each iteration is $O(m^3 (n+1)^3)$.} 

  \subsection{Urban pollution monitoring}\label{sec:sec:poll}
\begin{figure}[htbp] 
\centering
  \includegraphics[width=0.95\columnwidth]{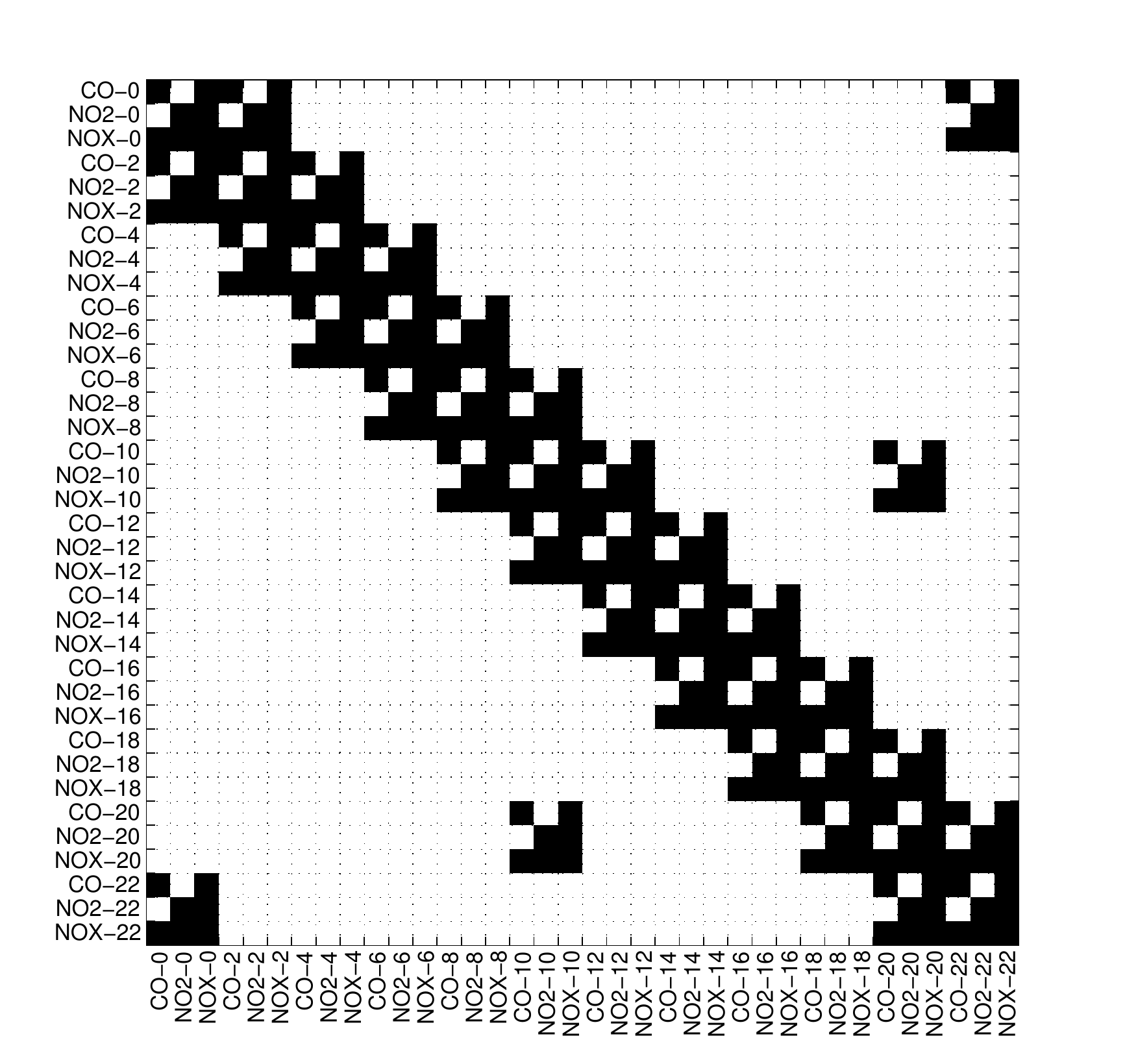}  
\caption{Sparsity pattern of the dynamic spatio-temporal graph for the urban atmospheric pollutants. Each node is denoted by A-{\mz hr} where A is the name of the pollutants and hr is the hour. A black square means that there is an edge between the corresponding nodes, otherwise a white square means there is not.} \label{Fig_SP}
\end{figure}

\begin{figure}[htbp] 
\centering
\subfloat[]{  \includegraphics[width=0.6\columnwidth]{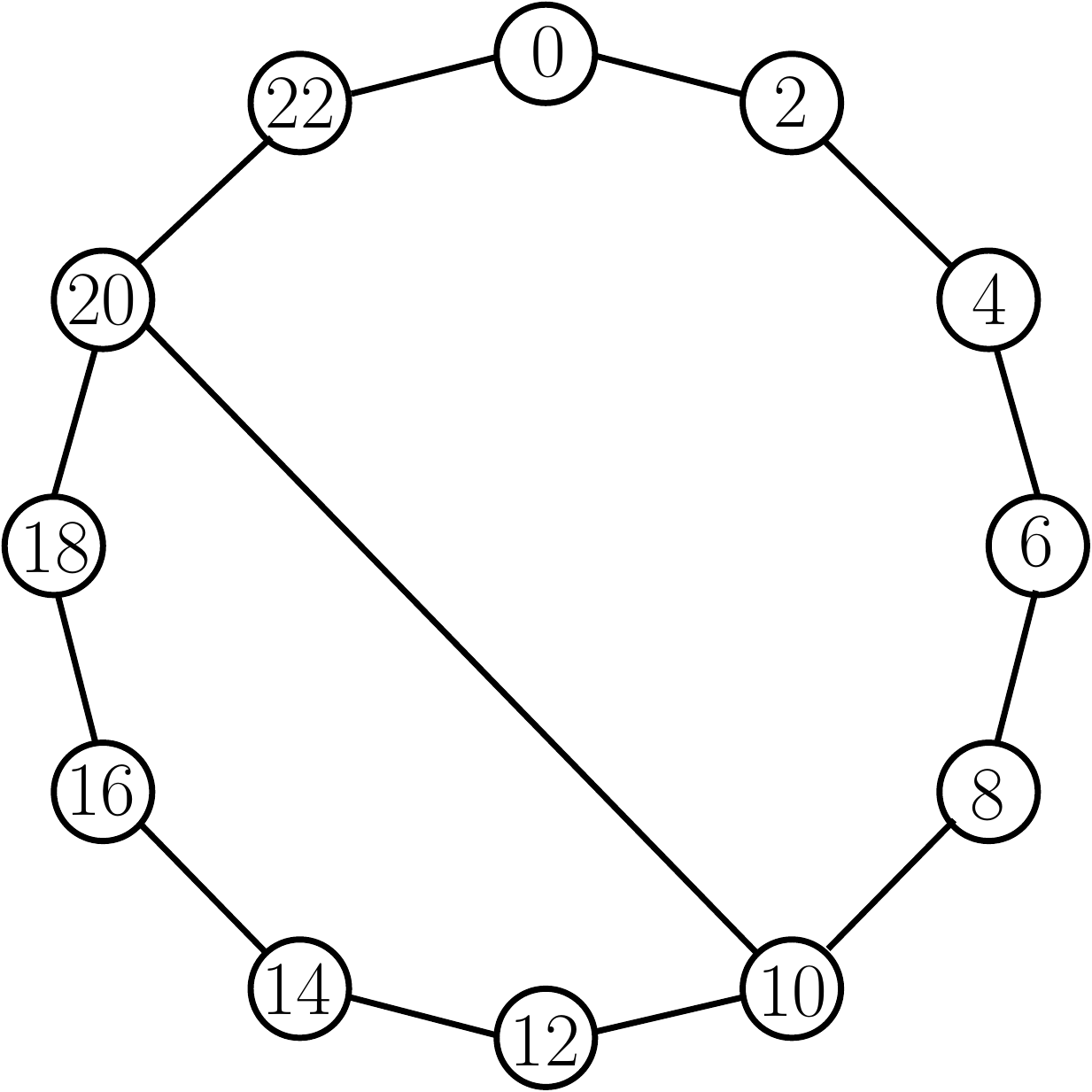}    }
\; \; \subfloat[]{ \includegraphics[width=0.3\columnwidth]{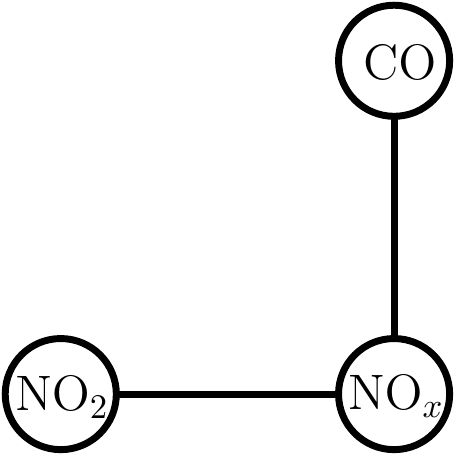} }
\caption{(a) Graphical model $\Gc(\Vc_1,\Ec_1)$ for the concentrations over a day (the number of the node corresponds to the hour). (b) Graphical model $\Gc(\Vc_2,\Ec_2)$ for the three urban atmospheric pollutants. } \label{Fig_SP2}
\end{figure}

\begin{figure}[htbp] 
\centering
  \includegraphics[width=0.95\columnwidth]{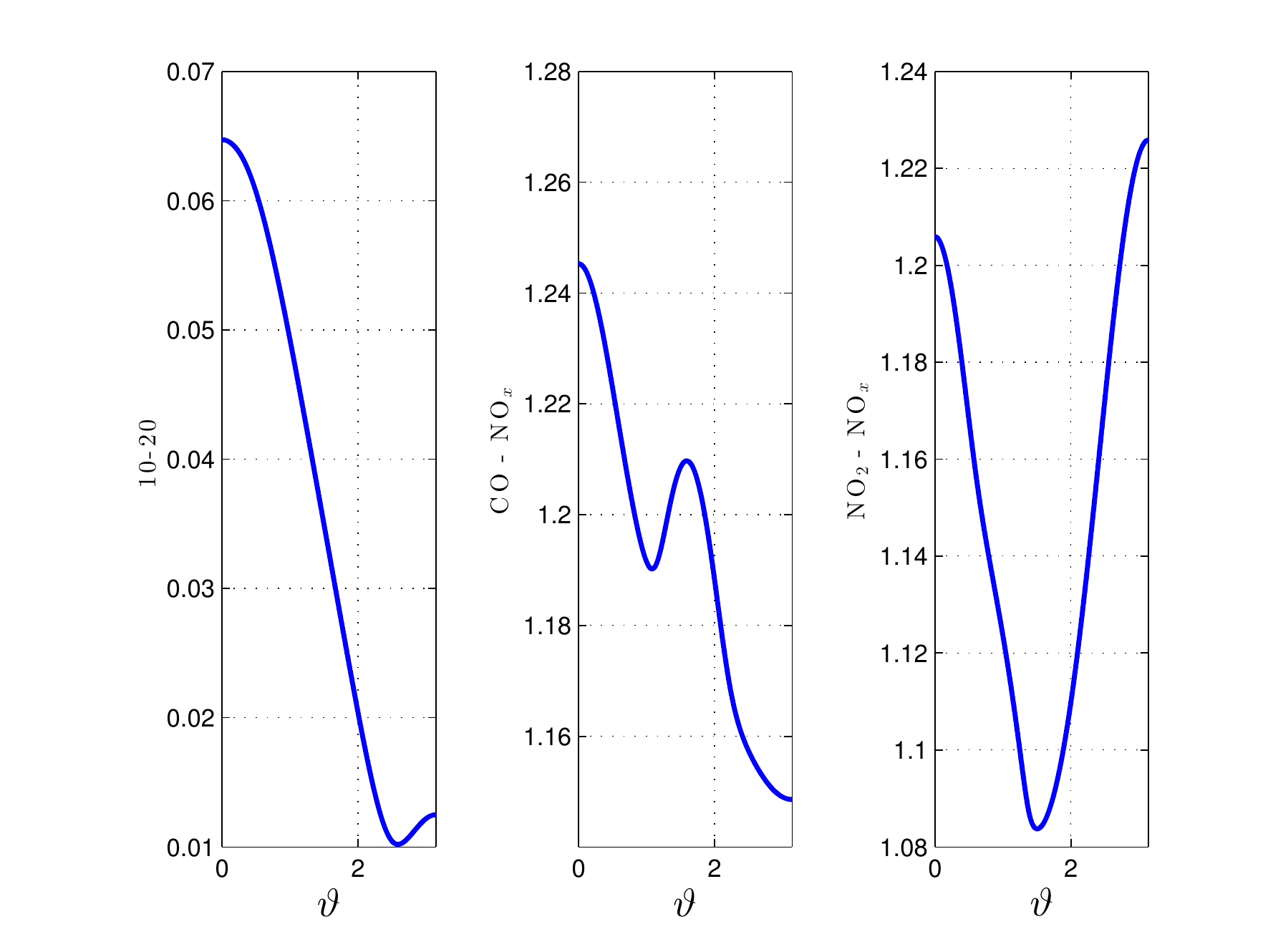}  
\caption{Norm of the spectra of the residuals characterizing the edge between the concentrations at 10 and at 20 (left); the edge between CO and NO$_x$ (center); the edge between NO$_2$ and NO$_x$ (right); } \label{Fig_partcor}
\end{figure}

We consider the concentration of the three pollutants CO, NO$_2$ and NO$_\mathrm x$ at a main street located in the center of an Italian city characterized by heavy car traffics. The corresponding three time series have been collected in the period 11 March 2004 - 3 April 2005 (389 days in total) by the regional environmental protection agency (ARPA) with sampling time equal to 1 hour, for more details see \cite{de2009co}. {\mz We aggregate the data in order to obtain time series of the averaged concentrations with sampling time 2 hours.} We normalize each time series in such a way that its sample variance is equal to one. These data describe the three-dimensional stochastic process $x=\{x(t),\; t\in\Zs \}$
with $x(t)=[\, x_1(t)\; x_2(t)\; x_3(t)\, ]^T$, $x_1(t)$ denotes the {\mz average} concentration  of CO at time $t$, $x_2(t)$ denotes the {\mz average} concentration  of NO$_2$ at time $t$ and {\mz $x_3(t)$ denotes the {\mz average} concentration  of NO$_\mathrm x$} at time $t$. $x(t)$ is non-stationary during a day: the peak time behaviour will be different to the off-peak behaviour. Therefore, we consider the process               
\al{y(t)=[\, x(12(t-1)+1)^T\ldots\, x(12 t )^T\, ]^T}
taking values in $\Rs^{36}$ and the corresponding sampling time is equal to 1 day. In this way we have a dataset $y^N$ with $N=389$. Once the dataset $y^N$  has been detrended, we apply method K1 with $m_1=12$, $m_2=3$ and the order of the AR process is set equal to $n=2$. The sparsity pattern of the estimated dynamic spatio-temporal graph is depicted in Figure \ref{Fig_SP}.  In view of Proposition \ref{prop_MEG}, we can characterize the graphical model for the three urban pollutants, $\Gc(\Vc_2,\Ec_2)$, and for their {\mz average} concentrations every 2 hours over a day, $\Gc(\Vc_1,\Ec_1)$, see Figure \ref{Fig_SP2}. 
 In regard to $\Gc(\Vc_1,\Ec_1)$, as expected, adjacent hours are conditionally dependent. Moreover, the concentrations at 10 and at 20 are conditionally dependent. The latter could explain the work journey (with a delay of 2-3 hours): people starts to work around 8 and finishes around 17. Finally, in Figure \ref{Fig_partcor}   we show the norm of the spectra $\tilde \Phi_\varepsilon $ (see the proof of Proposition \ref{prop_MEG}) characterizing the edge connecting the concentrations at 10 and at 20 in  $\Gc(\Vc_1,\Ec_1)$, and the two edges in  $\Gc(\Vc_2,\Ec_2)$. For all these edges, most of conditional dependence happens at low frequencies. We also applied method K2: we have obtained the same graph topology.

\begin{figure}[htbp] 
\centering
\subfloat[]{  \includegraphics[width=0.6\columnwidth]{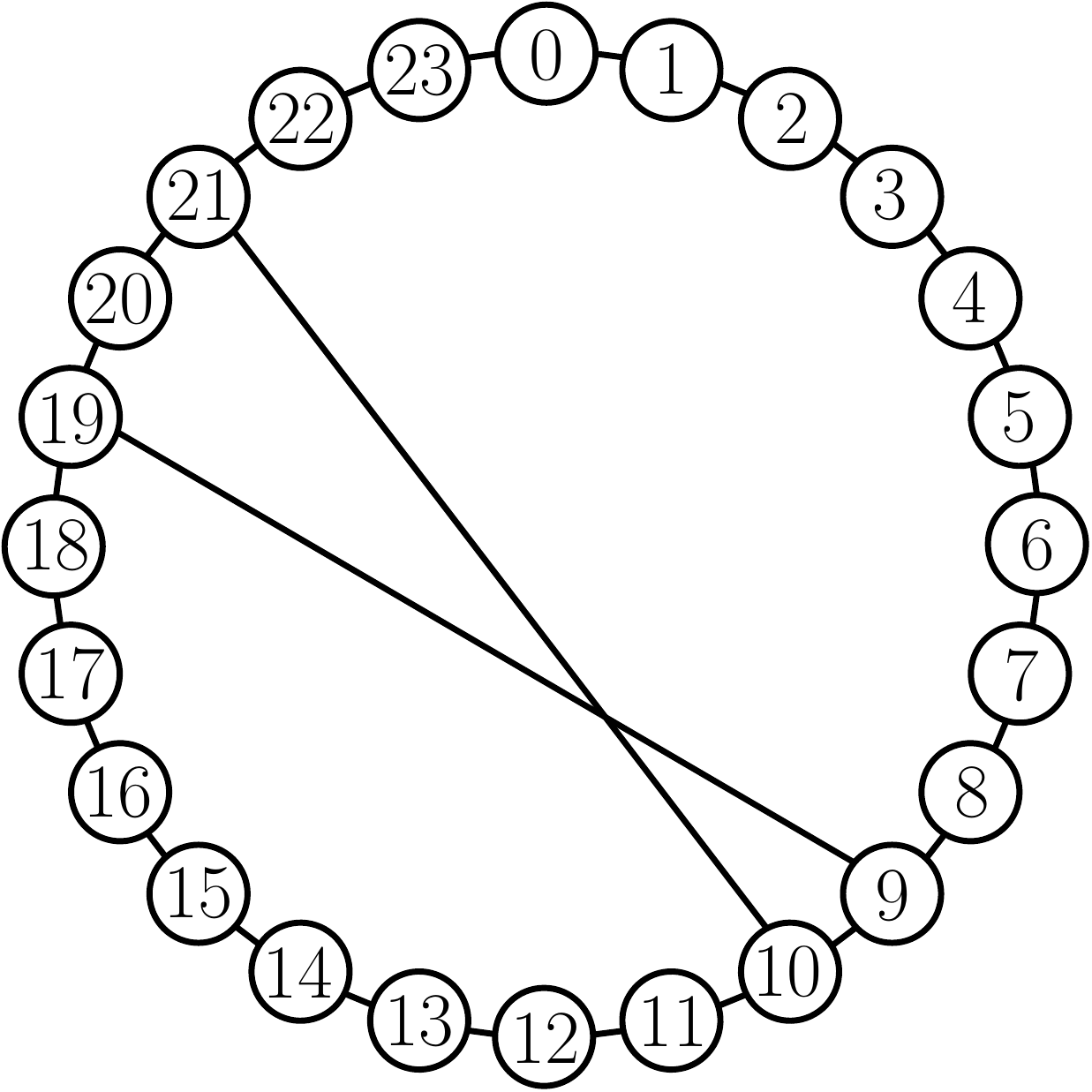}    }
\; \; \subfloat[]{ \includegraphics[width=0.3\columnwidth]{e1.pdf} }
\caption{(a) Graphical model $\Gc(\Vc_1,\Ec_1)$ for the concentrations over a day (the number of the node corresponds to the hour). (b) Graphical model $\Gc(\Vc_2,\Ec_2)$ for the three urban atmospheric pollutants. } \label{Fig_SP224}
\end{figure}

{\mz Finally, as sanity check, we estimate the KGM using the original data, i.e. the ones with sampling time equal to 1 hour. We apply method K1 with $m_1=24$, $m_2=3$ and $n=4$. It is worth noting that the order of this AR process has been chosen in such a way that its value at a certain time depends on its past values over a time interval of 4 hours as in the one considered for the aggregated data. The topology of the estimated KGM is depicted in Figure \ref{Fig_SP224}. Such a model is consistent with the one obtained by the aggregated data: the graphical model for the three urban pollutants is the same; adjacent hours are conditionally dependent; there are some conditional dependence relations between the hours 9-10 and 19-21. }

\section{Conclusions} \label{sec:concl}
 We have introduced a KGM corresponding to an AR Gaussian stochastic process. The latter is described by a PSD whose inverse has support which can be decomposed as a Kronecker product. We have proposed a ML estimator for KGM adopting a Bayesian perspective. In particular we have introduced two priors for the estimation of the hyperparameters: the max prior of Section \ref{sec:additive} and the multiplicative prior of Section \ref{sec:additive2}. Although the latter has been successfully used for collaborative filtering and multi-task learning, it provides a performance which is worse than the one of the max prior. We have also shown that the ML estimator is connected to a ME problem. Finally, we have tested the proposed approach to synthetic data as well as urban pollution data.

\end{document}